\documentclass[10pt]{article}
\usepackage{amssymb,amsmath,amsthm}
\usepackage{indentfirst}
\usepackage{exscale}
\usepackage{relsize}

\theoremstyle{definition}

\theoremstyle{remark}

\numberwithin{equation}{section}
\begin{document}

\title{\Large\bf{Existence and multiplicity of solutions for a class of  Kirchhoff type $(\Phi_1,\Phi_2)$-Laplacian system with locally super-linear condition in $\mathbb{R}^N$}
 }
\date{}
\author {\ Cuiling Liu, \ Xingyong Zhang\footnote{Corresponding author, E-mail address: zhangxingyong1@163.com}\\
      {\footnotesize Faculty of Science, Kunming University of Science and Technology, Kunming, Yunnan, 650500, P.R. China.}\\
      }
 \date{}
 \maketitle

 \begin{center}
 \begin{minipage}{15cm}

\small  {\bf Abstract:}
We investigate  the existence and multiplicity of weak solutions for a nonlinear Kirchhoff type quasilinear elliptic system on the whole space $\mathbb{R}^N$. We assume that the nonlinear term satisfies the locally super-$(m_1,m_2)$ condition, that is,
 $\lim_{|(u,v)|\rightarrow+\infty}\frac{F(x,u,v)}{|u|^{m_1}+|v|^{m_2}}=+\infty \mbox{\ \ for a.e.\ \ } x \in G$
 where $G$ is a domain in $\mathbb{R}^N$, which is weaker than the well-known Ambrosseti-Rabinowitz condition and the naturally global restriction, $\lim_{|(u,v)|\rightarrow+\infty}\frac{F(x,u,v)}{|u|^{m_1}+|v|^{m_2}}=+\infty \mbox{\ \ for a.e.\ \ } x \in \mathbb{R}^N$.  We obtain that system has at least one weak solution by using the classical Mountain Pass Theorem. To a certain extent, our theorems extend the results of Tang-Lin-Yu [Journal of Dynamics and Differential Equations, 2019, 31(1): 369-383]. Moreover, under the above naturally global restriction, we obtain that system has infinitely many weak solutions of high energy by using the Symmetric Mountain Pass Theorem, which is different from those results of Wang-Zhang-Fang [Journal of Nonlinear Sciences and Applications, 2017, 10(7): 3792-3814] even if we consider the system on the bounded domain with Dirichlet boundary condition.

\end{minipage}
 \end{center}
  \allowdisplaybreaks
\section{Introduction}
In this paper, we are dedicated to studying the existence and multiplicity of weak solutions for the following generalized nonlinear and non-homogeneous Kirchhoff type elliptic system in Orlicz-Sobolev spaces
 \begin{equation}\label{eq1}
 \left\{
  \begin{array}{ll}
 -M_{1}\left(\int_{\mathbb{R}^{N}}\Phi_1(|\nabla u|)dx\right)\Delta_{\Phi_1}u+V_{1}(x)\phi_{1}(|u|)u=F_u(x,u,v), &x\in \mathbb{R}^N,\\
 -M_{2}\left(\int_{\mathbb{R}^{N}}\Phi_2(|\nabla v|)dx\right)\Delta_{\Phi_2}v+V_{2}(x)\phi_{2}(|v|)v=F_v(x,u,v), &x\in \mathbb{R}^N,\\
  u\in W^{1,\Phi_{1}}(\mathbb{R}^N),v\in W^{1,\Phi_{2}}(\mathbb{R}^N),
    \end{array}
 \right.
 \end{equation}
where $\Delta_{\Phi_i}(u)=\mbox{div}(\phi_i(|\nabla u|)\nabla u)$, $(i=1, 2)$, $\phi_i:(0,+\infty)\rightarrow(0,+\infty)$ are two functions
which satisfy the following conditions:
\begin{itemize}
	 \item[$(\phi_1)$] $\phi_i\in C^1[(0,+\infty),\mathbb{R}^+]$, $t\phi_i(t)\rightarrow0$ as
 $t\rightarrow0$, $t\phi_i(t)\rightarrow+\infty$ as
 $t\rightarrow+\infty$;
 \item[$(\phi_2)$] $t\rightarrow t\phi_i(t)$ are strictly increasing;
 \item[$(\phi_3)$] $1<l_i:=\inf_{t>0}\frac{t^2\phi_i(t)}{\Phi_i(t)}\le\sup_{t>0}\frac{t^2\phi_i(t)}{\Phi_i(t)}=:m_i<\min\{N,l_{i}^{\ast}\}$, where $\Phi_i(t):=\int_{0}^{|t|}s\phi_i(s)ds, \ t\in \mathbb{R}$ and $l_{i}^{\ast}=\frac{l_{i}N}{N-l_{i}}$;
 \item[$(\phi_4)$]
    there exist positive constants $c_{i,1}$ and $c_{i,2}$, $i=1,2$, such that
    $$
    c_{i,1}|t|^{l_i}\le \Phi_i(t)\le c_{i,2}|t|^{l_i},\ \ \forall |t|<1;
    $$
\end{itemize}
\par
 Moreover, we introduce the following conditions on $F$, $V_{i}$ and $M_{i}$:
\begin{itemize}
\item[$(F_0)$]  $F: \mathbb{R}^{N}\times \mathbb{R}\times \mathbb{R}\rightarrow \mathbb{R}$ is a $C^1$ function such that
                $F(x,0,0)=0$ for all $x\in \mathbb{R}^{N}$ and $F(x,u,v)\ge 0$ for all $(x,u,v)\in \mathbb{R}^N\times \mathbb{R}\times \mathbb{R}$;
\end{itemize}
\begin{itemize}
\item[$(V_0)$]  $V_{i} \in C(\mathbb{R}^{N},\mathbb{R})$ and $\inf_{\mathbb{R}^{N}}V_{i}(x)> 1$,\; $i=1,2$;
\end{itemize}
\begin{itemize}
\item[$(V_1)$]  there exist  constants $C_{i,1}>0$ such that
                $$
                   \lim_{|z|\rightarrow \infty} \mbox{meas}\{x \in \mathbb{R}^{N}: |x-z|\leq C_{i,1},
                          V_{i}(x)\leq C_{i,2}\}=0\;\;\mbox{for every}\; C_{i,2}>0, i=1,2,
                $$
                where meas$(\cdot)$ denotes the Lebesgue measure in $\mathbb{R}^{N}$;
\end{itemize}
\begin{itemize}
\item[$(M_0)$]  $M_{i} \in C(\mathbb{R}^+,\mathbb{R}^+) $  and $C_{i,3} \leq M_{i}(t)\leq C_{i,4},\; \forall \;t\geq 0$ for some $C_{i,3},C_{i,4}>0$,\; $i=1,2$;
\item[$(M_1)$]  $\widehat{M_{i}}(t):=\int_{0}^{t}M_{i}(s)ds\geq M_{i}(t)t,\; i=1,2$.
\end{itemize}
\par
Let $\phi_1=\phi_2=:\phi$, $v=u$,  $M_1=M_2=:M$, $V_1=V_2=:V$ and $F(x,u,v)=F(x,v,u)$.
Then the system \eqref{eq1} reduces to the following non-homogeneous and nonlocal quasilinear elliptic equation
 \begin{equation}\label{eq1.2}
 \left\{
  \begin{array}{ll}
 -M\left(\int_{\Omega}\Phi(|\nabla u|)dx\right)\Delta_{\Phi}u+V(x)\phi(|u|)u=f(x,u), &x\in  \Omega,\\
  u\in W^{1,\Phi}(\Omega),
    \end{array}
 \right.
 \end{equation}
where $\Omega$ is a domain in $\mathbb{R}^N$ and $f(x,u)=F_u(x,u,u)$.
\par
In recent years, many authors are concerned with nonlocal problems like \eqref{eq1.2} which can been seen as a generalization of the second-order semilinear elliptic equations, $p$-Laplacian equations and $(p,q)$-Laplacian equations (see \cite{Tang2016,Alves2005,Claudianor2012,ZhangQF2021,ZhangQF20212,Sbai2021,Tsouli2019,Nguyen2013,Figueiredo2015, Heidari2021,Ferrara2016,Chung2014,Adriouch2006,Alves2014, Kim2020,Lee2020,Santos2019,Avci2020, Carvalho2015, Liu-Shibo2019}
and references therein).

\par
 Next, we emphasize the results in \cite{Tang2019} which mainly inspires our works in this paper. In \cite{Tang2019}, Tang,  Lin and Yu  investigated the existence of nontrivial solutions for the following semilinear Schr$\ddot{\mbox{o}}$dinger equation:
\begin{equation}\label{eq1.4}
 \left\{
  \begin{array}{ll}
 -\Delta u+V(x)u=f(x,u),\; x\in \mathbb{R}^N,\\
 u\in H^{1}(\mathbb{R}^N),
    \end{array}
 \right.
 \end{equation}
where the potential $V\in C(\mathbb{R}^N,\mathbb{R})$ is a sign-changing function which satisfies the periodic or coercive conditions. If
$f$  satisfies a subcritical condition and  the following  locally super-quadratic condition:
\par
\noindent
{\it $(f_1)$ there exists a domain $G\subset \mathbb{R}^N$ such that
$$
\lim_{t\rightarrow\infty}\frac{\int_{0}^{t}f(x,s)ds}{t^{2}}=+\infty, \; a.e.\; x\in G,
$$
} which is weaker than the following naturally super-quadratic condition:
\par
\noindent
$(f'_1)\;\;\lim_{t\rightarrow\infty}\frac{F(x,t)}{|t|^{2}}=+\infty \; \mbox{\it  uniformly for}\; x\in \mathbb{R}^N,$
\par
\noindent
by using the linking geometry theorem,  they obtained that the problem \eqref{eq1.4} has a nontrivial weak solution.
\par
There have been some contributions devoted to the study of system (\ref{eq1}) with $M_i=1$ involving the existence and multiplicity of weak solutions.
In \cite{wang2017}, Wang, Zhang  and Fang considered the following quasilinear elliptic system in Orlicz-Sobolev spaces:
\begin{equation}\label{eq1.3}
 \left\{
  \begin{array}{ll}
 -\mbox{div}(\phi_1(|\nabla{u}|)\nabla{u})=F_u(x,u,v), &\mbox{in } \Omega,\\
 -\mbox{div}(\phi_2(|\nabla{v}|)\nabla{v})=F_v(x,u,v), &\mbox{in } \Omega,\\
 u=v=0, &\mbox{on } \partial\Omega,
    \end{array}
 \right.
 \end{equation}
where $\Omega$ is a bounded domain in $\mathbb{R}^N(N\geq2)$ with smooth
boundary $\partial\Omega$. When $F$ satisfies some appropriate conditions including $(\phi_1,\phi_2)$-superlinear and subcritical growth conditions at infinity  as well as symmetric condition, by using the mountain pass theorem and the symmetric mountain pass theorem, they obtained that system \eqref{eq1.3} has a nontrivial weak solution and infinitely many weak solutions,  respectively. Subsequently, more works were obtained for systems like  (\ref{eq1}) (see \cite{wang2017-2}-\cite{Zhang2022}).
For example, in \cite{wang2017-4}, Wang, Zhang and Fang considered the  quasilinear elliptic system like \eqref{eq1} with $M_{i}(t)=1$ in $\mathbb{R}^{N}$.
When the potential functions are bounded, $F$ satisfies sub-linear growth  condition, by using the least action principle, they obtained that system has at least one nontrivial weak solution. If  $F$ also satisfies a symmetric condition, by using the genus theory, they obtained that system has infinitely many weak solutions. In \cite{Zhang2021}-\cite{Zhang2022}, we developed the Moser iteration technique, and then by using the mountain pass theorem and cut-off technique, we obtain that system like \eqref{eq1} with $M_{i}(t)=1$ and a parameter $\lambda$ has a nontrivial weak solution $(u_\lambda,v_\lambda)$ with  $\|(u_{\lambda},v_{\lambda})\|_{\infty}\le 2$  for every $\lambda$ large enough if the nonlinear term $F$ satisfies some growth conditions only in a circle with center $0$ and radius $4$.

\par
Inspired by  \cite{Alves2014,Carvalho2015,Tang2019,wang2017}, in this paper, we shall investigate the existence and multiplicity of weak solutions for system \eqref{eq1} with locally super-$(m_1,m_2)$ growth in $\mathbb{R}^{N}$.
We assume that $F$ satisfies the following local condition
 $$
 \lim_{|(u,v)|\rightarrow+\infty}\frac{F(x,u,v)}{|u|^{m_1}+|v|^{m_2}}=+\infty \mbox{\ \ for a.e.\ \ } x \in G,
 $$
where $G$ is a domain in $\mathbb{R}^{N}$, by using  the Mountain Pass Theorem,  we obtain that system \eqref{eq1} has a nontrivial weak solution. Besides,
if $F$ also satisfies a symmetric condition and the following naturally global restriction
 $$
 \lim_{|(u,v)|\rightarrow+\infty}\frac{F(x,u,v)}{|u|^{m_1}+|v|^{m_2}}=+\infty \mbox{\ \ for a.e.\ \ } x \in \mathbb R^N,
 $$
 by using the Symmetric Mountain Pass Theorem, we obtain that system \eqref{eq1} has infinitely many weak solutions of high energy.
We develop some results in some known references in the following sense:
\begin{itemize}
\item[(\uppercase\expandafter{\romannumeral1})]
Our local condition extends the locally super-quadratic condition  of \cite{Tang2019} (see $(f_1)$ mentioned above);
\item[(\uppercase\expandafter{\romannumeral2})] Different from those in  \cite{wang2017}-\cite{Zhang2022}, we  consider the nonlocal Kirchhoff-type problems;
\item[(\uppercase\expandafter{\romannumeral3})]
Our conditions are weaker than the Ambrosetti-Rabinowitz ((A-R) for short) condition in \cite{Alves2014,Liu-Shibo2019};
\item[(\uppercase\expandafter{\romannumeral4})]
 Different from that in \cite{wang2017}, we work in the whole-space $\mathbb{R}^N$ rather than a bounded domain $\Omega\subset\mathbb{R}^N$. Especially, we introduce a new $(\phi_1,\phi_2)$-superlinear condition (see $(F_4)$ below and  Remark 3.4 for details), which is different from the following condition in \cite{wang2017} even if we restrict $(F_4)$ to the bounded domain $\Omega$:\\
{\it $(f_2)$ there exists a continuous function $\gamma: [0,\infty)\rightarrow \mathbb{R}$ and it satisfies that $\Gamma(t):=\int_{0}^{|t|}\gamma(s)ds,~ t\in\mathbb{R}$ is an $N$-function with
 \begin{equation*}
 1<l_\Gamma:=\inf_{t>0}\frac{t\gamma(t)}{\Gamma(t)}\leq \sup_{t>0}\frac{t\gamma(t)}{\Gamma(t)}=:m_\Gamma<+\infty,
 \end{equation*}
 such that
 \begin{equation*}\label{3.1.3+}
 \Gamma\left(\frac{F(x,u,v)}{|u|^{l_1}+|v|^{l_2}}\right)\leq d_{1}\overline{F}(x,u,v), \quad x\in\Omega, \quad |(u,v)|\geq r_{1},
 \end{equation*}
 where constants $d_{1}, r_{1}>0$ and
 $$\overline{F}(x,u,v):=\frac{1}{m_1}F_u(x,u,v)u+\frac{1}{m_2}F_v(x,u,v)v-F(x,u,v), ~~~\forall  (x,u,v)\in \Omega\times \mathbb{R}\times \mathbb{R},$$
and  the following $(\phi_1,\phi_2)$-superlinear growth conditions hold:
$$
 \lim_{|(u,v)|\rightarrow+\infty}\frac{F(x,u,v)}{\Phi_1(u)+\Phi_2(v)}=+\infty \mbox{\ \ uniformly for all } x \in \Omega,
$$
where $\Omega$ is a bounded domain in $\mathbb{R}^N$;}
\item[(\uppercase\expandafter{\romannumeral5})]
Similar to (\uppercase\expandafter{\romannumeral4}), for the scalar equation,  we  also obtain  a new $\phi$-superlinear condition which is different from those in  \cite{Carvalho2015,wang2017} even if we restrict it to the bounded domain $\Omega$. One can see  the Remark  4.3
 for details;
\item[(\uppercase\expandafter{\romannumeral6})]
Because of the coupling relationship of $u$ and $v$ and the inhomogeneous properties of $\Phi_i,i=1,2$, our proofs become more difficult and complex than those in \cite{Tang2019}. Especially, such difficulty and complexity can be embodied in the proofs of the compactness of Cerami sequence. Moreover, because the new condition $(F_4)$ below is different from  $(f_2)$ and we consider the problem (\ref{eq1}) on the whole space $\mathbb{R}^N$ rather than a bounded domain $\Omega$, our proofs on the compactness of Cerami sequence are different from those in \cite{wang2017}.
\end{itemize}
\par
The remainder of this article focuses on some preliminaries, the main results of this paper and their proofs and an example which illustrates our results. Finally, a remark on semi-trivial solutions of  (\ref{eq1}) is given.

\section{Preliminaries}\label{section 2}
In this section, to deal with such problem for system \eqref{eq1}, we need to briefly list some fundamental definitions and essential properties of Orlicz and Orlicz-Sobolev spaces   and introduce some classical results from variational methods.  For a deeper understanding of these concepts, we refer readers for more details to the books \cite{Adams2003,Rabinowitz1986,M. M. Rao2002}.
\vskip2mm
\noindent
{\bf Definition 2.1.} \cite{Adams2003} Let $b$ : $[0,+\infty)\rightarrow [0,+\infty)$ be a right continuous, monotone increasing function with
\begin{itemize}
	 \item[$\rm(1)$] $b(0)=0$;
 \item[$\rm(2)$] $\lim_{t\rightarrow +\infty} b(t)=+\infty$;
 \item[$\rm(3)$] $b(t)>0 $ whenever $t>0$.
\end{itemize}
 Then the function defined on $\mathbb{R}$ by $B(t)=\int_{0}^{|t|}b(s)ds$
 is called as an $N$-function.
 \vskip2mm
 \par
  By the definition of $N$-function $B$, it is obvious that $B(0)=0$ and $B$ is strictly convex. We call that an $N$-function $B$ satisfies a $\Delta_2$-condition globally (or near infinity) if
 $$
  \sup_{t>0}\frac{B(2t)}{B(t)}<+\infty \ \ (\mbox{ or } \limsup_{t\rightarrow \infty}\frac{B(2t)}{B(t)}<+\infty),
 $$
  which implies that there exists a constant $K>0$ such that $B(2t)\leq K B(t)$ for all $t\geq0$ (or $t\geq t_0>0$). We also state the equivalent form that $B$ satisfies a $\Delta_2$-condition globally (or near infinity) if and only if for any $c\geq 1$, there exists a constant $K_c>0$ such that $B(ct)\leq K_c B(t)$ for all $t\geq0$ (or $t\geq t_0>0$).
  \vskip2mm
  \noindent
{\bf Definition 2.2.}\cite{Adams2003}  For an $N$-function $B$, we define $$ \widetilde{B}(t)=\int_{0}^{|t|}b^{-1}(s)ds, \quad t\in \mathbb{R},$$
 where $b^{-1}$ is the right inverse of the right derivative $b$ of $B$. Then $\widetilde{B}$ is an $N$-function
called as the complement of $B$.
\vskip2mm
\par
 It holds that Young's inequality (see \cite{Adams2003,M. M. Rao2002})
\begin{equation}\label{2.1.1}
st\leq B (s)+\widetilde{B}(t), \quad s, t\geq 0
\end{equation}
and the inequality (see \cite[Lemma A.2]{Fukagai2006})
 \begin{equation}\label{2.1.2}
 \widetilde{B}(b(t))\leq B(2t), \quad t\geq 0.
 \end{equation}
\par
Now, we recall the Orlicz space $L^{B}(\Omega)$ associated with $B$.
The Orlicz space $L^{B}(\Omega)$ is the vectorial space of the measurable functions $u: \Omega\rightarrow \mathbb{R}$ satisfying
$$\int_{\Omega}B(|u|)dx<+\infty,$$
where $\Omega \subset \mathbb{R}^N$ is an open set. $L^{B}(\Omega)$ is a Banach space endowed with Luxemburg norm
$$
\|u\|_{B}:=\inf \left\{\lambda >0: \int_{\Omega}B \left(\frac{u}{\lambda}\right)dx\le 1\right\}.
$$
 The fact that $B$ satisfies $\Delta_2$-condition globally implies that
\begin{equation}\label{2.1.3+}
 u_n\rightarrow u \mbox{ in } L^{B}(\Omega) \Longleftrightarrow \int_{\Omega}B(u_n-u)dx\rightarrow 0.
 \end{equation}
 Moreover, a generalized type of H\"{o}lder's inequality (see \cite{Adams2003,M. M. Rao2002})
$$\left| \int_{\Omega}uvdx \right|\leq 2\|u\|_{\Phi}\|v\|_{\widetilde{\Phi}}, \quad \mbox{ for all } u \in L^{\Phi}(\Omega)\mbox{ and } v \in L^{\widetilde{\Phi}}(\Omega)$$
can be gained by applying Young's inequality \eqref{2.1.1}.
\par
The corresponding Orlicz-Sobolev space (see \cite{Adams2003,M. M. Rao2002}) is defined by
$$
W^{1, B}(\Omega):=\left\{u \in L^{B}(\Omega): \frac{\partial u}{\partial x_i} \in L^{B}(\Omega), i=1,\cdots, N\right\}
$$
with the norm
$$\|u\|_{1, B}:=\|u\|_{B}+\|\nabla u \|_{B}.$$
\par
Consider the subspace $X$ of $W^{1,B}(\mathbb{R}^N)$,
     \begin{eqnarray}\label{dddc1}
  X=\left\{u\in W^{1,B}(\mathbb{R}^N)\Big| \int_{\mathbb{R}^N}V(x)B(|u|)dx<\infty\right\}
  \end{eqnarray}
  with the norm
  \begin{eqnarray}\label{2.1.4}
  \|u\|_X=\|\nabla u\|_B+\|u\|_{B,V},
  \end{eqnarray}
  where
  $$
  \|u\|_{B,V}=\inf\left\{\alpha>0\Big|\int_{\mathbb{R}^N} V(x)B\left(\frac{|u|}{\alpha}\right)dx\le 1\right\}
  $$
  and $\inf_{x\in\mathbb{R}^N} V(x)>0$. Then $(X,\|\cdot\|)$ is a separable and reflexive Banach space (see \cite{Liu-Shibo2019}).
\vskip2mm
\noindent
{\bf Lemma 2.3.}\cite{Adams2003,Fukagai2006}
  If $B$ is an $N$-function, then the following conditions are equivalent:
\begin{itemize}
	\item[$\rm(1)$]
\begin{equation}\label{2.1.5}
1\leq l=\inf_{t>0}\frac{tb(t)}{B(t)}\leq\sup_{t>0}\frac{tb(t)}{B(t)}=m<+\infty;
\end{equation}
\item[$\rm(2)$] Let $\zeta_0(t)=\min\{t^l, t^m\}$ and $\zeta_1(t)=\max\{t^l, t^m\}$, $t\geq0$. $B$ satisfies
$$\zeta_0(t)B(\rho)\leq B(\rho t)\leq \zeta_1(t)B(\rho), \quad \forall \rho, t\geq 0;$$
\item[$\rm(3)$] $B$ satisfies a $\Delta_2$-condition globally.
\end{itemize}
\noindent
{\bf Lemma 2.4.}\cite{Fukagai2006}
 If $B$ is an $N$-function and \eqref{2.1.5} holds, then $B$ satisfies
$$\zeta_0(\|u\|_{B})\leq\int_{\Omega}B(u)dx\leq\zeta_1(\|u\|_{B}), \quad \forall u \in L^{B}(\Omega).$$
\noindent
{\bf Lemma 2.5.}\cite{Fukagai2006}
 If $B$ is an $N$-function and \eqref{2.1.5} holds with $l>1$. Let $\widetilde{B}$ be the complement of $B$ and $\zeta_2(t)=\min\{t^{\widetilde{l}},t^{\widetilde{m}}\}$, $\zeta_3(t)=\max\{t^{\widetilde{l}},t^{\widetilde{m}}\}$ for $t\geq0$, where $\widetilde{l}:=\frac{l}{l-1}$ and $\widetilde{m}:=\frac{m}{m-1}$. Then $\widetilde{B}$ satisfies
\begin{itemize}
	\item[$\rm(1)$]
$$\widetilde{m}=\inf_{t>0}\frac{t\widetilde{B}^{'}(t)}{\widetilde{B}(t)}\leq\sup_{t>0}\frac{t\widetilde{B}^{'}(t)}{\widetilde{B}(t)}=\widetilde{l};$$
\item[$\rm(2)$]
$$\zeta_2(t)\widetilde{B}(\rho)\leq \widetilde{B}(\rho t)\leq \zeta_3(t)\widetilde{B}(\rho), \quad \forall  \rho, t\geq 0;$$
\item[$\rm(3)$]
$$\zeta_2(\|u\|_{\widetilde{B}})\leq\int_{\Omega}\widetilde{B}(u)dx\leq\zeta_3(\|u\|_{\widetilde{B}}), \quad \forall u \in L^{\widetilde{B}}(\Omega).$$
\end{itemize}
{\bf Lemma 2.6.}\cite{Fukagai2006}
 If $B$ is an $N$-function and \eqref{2.1.5} holds with $l,m\in(1, N)$. Let $\zeta_4(t)=\min\{t^{l^*},t^{m^*}\}$,\\ $\zeta_5(t)=\max\{t^{l^*},t^{m^*}\}$ for $t\geq0$, where $l^*:=\frac{lN}{N-l}$ and $m^*:=\frac{mN}{N-m}$. Then $B_{*}$ satisfies
\begin{itemize}
	\item[$\rm(1)$]
$$l^*=\inf_{t>0}\frac{t B_{*}'(t)}{B_*(t)}\leq\sup_{t>0}\frac{tB_{*}'(t)}{B_*(t)}=m^*;$$
\item[$\rm(2)$]
$$\zeta_4(t)B_*(\rho)\leq B_*(\rho t)\leq \zeta_5(t)B_*(\rho), \quad \forall \rho, t\geq 0;$$
\item[$\rm(3)$]
$$\zeta_4(\|u\|_{B_*})\leq\int_{\Omega}B_*(u)dx\leq\zeta_5(\|u\|_{B_*}), \quad \forall u \in L^{B_*}(\Omega),$$
where $B_*$ is the Sobolev conjugate function of $B$, which is defined by
$$B_{*}^{-1}(t)=\int_{0}^{t}\frac{B^{-1}(s)}{s^{\frac{N+1}{N}}}ds~~ \mbox{ \it for } t\geq 0~~~ \mbox{ \it and }~~~ B_*(t)=B_*(-t)~~ \mbox{ for } t\leq 0.$$
\end{itemize}

Next, we recall some embeddings.
Let $\Psi$ be an $N$-function verifying $\Delta_2$-condition. If
 \begin{equation}\label{2.1.6}
\overline{\lim_{t\rightarrow 0}}\frac{\Psi(t)}{B(t)}< +\infty\;\;\;\mbox{and}\;\;\;
\mathop{\overline{\lim}}_{|t|\rightarrow +\infty}\frac{\Psi(t)}{B_{\ast}(t)}< +\infty,
\end{equation}\label{3.1}
then we have a continuous embedding $W^{1,B}(\mathbb{R}^{N})\hookrightarrow L^{\Psi}(\mathbb{R}^{N})$.
Moreover, if
 \begin{equation}\label{2.1.7}
\lim_{|t|\rightarrow 0}\frac{\Psi(t)}{B(t)}< +\infty\;\;\;\mbox{and}\;\;\;
\lim_{|t|\rightarrow +\infty}\frac{\Psi(t)}{B_{\ast}(t)}=0,
 \end{equation}
then the embedding $W^{1,B}(\mathbb{R}^{N})\hookrightarrow L_{loc}^{\Psi}(\mathbb{R}^{N})$ is compact and we call that such $\Psi$ satisfies the subcritical condition.

\vskip2mm
\noindent
{\bf Lemma 2.7.}\cite{Liu-Shibo2019}
Assume that $b:[0,\infty)\rightarrow [0,\infty) \in C^{1}$ and $V$ satisfies the following conditions:
\begin{itemize}
\item[$(a_0)$]  the function $t\rightarrow b(t)t$ is increasing in $(0,\infty)$;
\item[$(a_1)$]  there exist $l, m\in (1,N)$ such that
                 $$
                 l\leq \frac{b(|t|)t^{2}}{B(t)}\leq m \;\; \mbox{for all}\;\; t\neq 0,
                 $$
                   where $l\leq m < l^{\ast}$ and $B(t)=\int_{0}^{|t|}b(s)sds$;
\end{itemize}
\begin{itemize}
    \item[$(V^{'}_0)$]  $V \in C(\mathbb{R}^{N},\mathbb{R})$ and $\inf_{x\in\mathbb{R}^{N}}V(x)= V_{0}>0$;
\end{itemize}
\begin{itemize}
\item[$(V^{'}_1)$]  for all $C_{0}>0$, $\mu(V^{-1}(-\infty,C_{0}])<\infty$,
                where $\mu$ is the Lebesgue measure in $\mathbb{R}^{N}$.
\end{itemize}
Then for any $N$-function $\Psi$ satisfying $\Delta_2$-condition
and \eqref{2.1.7}, the embedding from $X$ into $L^{\Psi}(\mathbb{R}^{N})$ is compact. Specifically, $X$ into $L^{B}(\mathbb{R}^{N})$ is compact,
where $X$ is defined by (\ref{dddc1}).

\vskip2mm
\noindent
{\bf Remark 2.8.}
   By Lemma 2.3 and Lemma 2.5, assumptions $(\phi_1)$--$(\phi_3)$ show that $\Phi_i ~(i=1, 2)$ and $\widetilde{\Phi}_i ~(i=1, 2)$ are $N$-functions satisfying $\Delta_2$-condition globally. Thus $L^{\Phi_i}(\mathbb{R}^N) (i=1, 2)$ and $W^{1, \Phi_i}(\mathbb{R}^N) (i=1, 2)$ are separable and reflexive Banach spaces (see \cite{Adams2003,M. M. Rao2002}).

\vskip2mm
 \par
 By the end of this section, we  recall  the mountain pass theorem
 (see \cite[Theorem 2.2]{Rabinowitz1986}) and  the symmetric mountain pass theorem (see \cite[Theorem 9.12]{Rabinowitz1986}) which will be used to prove Theorem 3.1 and
 Theorem 3.9 in Section 3, respectively.
 \par
 We first recall that  $I \in C^1(E,\mathbb{R})$ satisfies the Palais-Smale condition ((PS)-condition for short) if any (PS)-sequence $\{u_n\}\subset E$ has a convergent subsequence, where (PS)-sequence $\{u_n\}$ means that
 \begin{equation*}
 I(u_n)\mbox{ is bounded}, \quad\|I'(u_n)\|\rightarrow 0, \mbox{ as } n\rightarrow \infty,
 \end{equation*}
  and we call that $I \in C^1(E,\mathbb{R})$ satisfies the Cerami-condition ((C)-condition for short) if any (C)-sequence $\{u_n\}\subset E$ has a convergent subsequence, where (C)-sequence $\{u_n\}$ means that
 \begin{equation}\label{2.2.1}
 I(u_n) \mbox{ is bounded  and } (1+\|u_n\|)\|I'(u_n)\|\rightarrow 0, \mbox{ as } n\rightarrow \infty.
 \end{equation}
 \par
 By  the discussion in \cite{Bartolo1983}, the (PS)-condition can be substituted with (C)-condition in the following Lemma 2.9 and Lemma 2.10.
\vskip2mm
\noindent
{\bf Lemma 2.9.}{\cite[Theorem 2.2]{Rabinowitz1986}}(Mountain Pass Theorem)
  Let $E$ be a real Banach space and $I \in C^1(E,\mathbb{R})$ satisfying {\rm(PS)}-condition. Suppose $I(0)=0$ and
\begin{itemize}
	 \item[$\rm(I_1)$] there are constants $\rho, \alpha$$ >0$ such that $I\mid _{\partial B_{\rho}}\geq \alpha$, and
 \item[$\rm(I_2)$] there is an $e \in E\backslash B_{\rho}$ such that $I(e)\leq0$.
\end{itemize}
 Then $I$ possesses a critical value $c\geq\alpha$.

\vskip2mm
\noindent
{\bf Lemma 2.10.}{\cite[Theorem 9.12]{Rabinowitz1986}} (Symmetric Mountain Pass Theorem)
    Let $E$ be an infinite-dimensional Banach space and let $I \in C^1(E,\mathbb{R})$ be even, satisfy {\rm(PS)}-condition, and $I(0)=0$. If $E=V\oplus X$, where $V$ is finite dimensional, and $I$ satisfies
\begin{itemize}
	 \item[$\rm(I_1)$] there are constants $\rho$, $\alpha$$ >0$ such that $I\mid _{\partial B_{\rho}\cap X}\geq \alpha$, and
 \item[$\rm(I_2)$] for each finite dimensional subspace $\widetilde{E}\subset E$, there is an $R=R(\widetilde{E})$ such that $I\leq 0$ on $\widetilde{E}\backslash B_{R(\widetilde{E})}$,
 where $B_R=\{u\in E: \|u\|<R\}$,
\end{itemize}
then $I$ possesses an unbounded sequence of critical values.

\section{Main results and proofs}
Define
 \begin{eqnarray}\label{ddd1}
 W_i=\left\{u\in W^{1,\Phi_i}(\mathbb{R}^N)\Big| \int_{\mathbb{R}^N}V(x)\Phi_i(|u|)dx<\infty\right\}
  \end{eqnarray}
  with the norm
  $$
      \|u\|_{i}
  =    \|\nabla u\|_{\Phi_i}+\|u\|_{\Phi_i,V_i}, \quad i=1,2.
  $$
Throughout this paper, we work in the subspace $W=:W_{1}\times W_{2}$ of $W^{1,\Phi_1}(\mathbb{R}^N)\times W^{1,\Phi_2}(\mathbb{R}^N)$ with the norm
  $$
        \|(u,v)\|
  =    \|u\|_{1}+ \|v\|_{2}
  =    \|\nabla u\|_{\Phi_1}+\|u\|_{\Phi_1,V_1}
     +\|\nabla v\|_{\Phi_2}+\|v\|_{\Phi_2,V_2}.
  $$
  Then $(W,\|\cdot\|)$ is a separable and reflexive Banach space.
\par
Define the energy functional $I$ on $W$ corresponding to system \eqref{eq1} is
 \begin{eqnarray}\label{3.0.1}
        I(u,v):
&= &       \widehat{M_{1}}\left(\int_{\mathbb{R}^{N}}\Phi_1(|\nabla u|)dx\right)
          + \widehat{M_{2}}\left(\int_{\mathbb{R}^{N}}\Phi_2(|\nabla v|)dx\right)
          + \int_{\mathbb{R}^{N}}V_{1}(x)\Phi_1(| u|)dx
          \nonumber\\
&&
          + \int_{\mathbb{R}^{N}}V_{2}(x)\Phi_2(| v|)dx
          - \int_{\mathbb{R}^{N}}F(x,u,v)dx,\quad (u,v)\in W.
 \end{eqnarray}
 Under the assumptions $(\phi_1)$--$(\phi_3)$, $(F_1)$, $(V_0)$, $(V_1)$, $(M_0)$ and $(M_1)$, by using the standard arguments as in \cite{wang2017,Huidobro1999}, we can prove that $I$ is well-defined and of class $C^1(W,\mathbb{R})$ with
 \begin{eqnarray}\label{3.0.2}
        \langle I'(u,v),(\tilde{u},\tilde{v})\rangle
 & = &      M_{1}\left(\int_{\mathbb{R}^{N}}\Phi_1(|\nabla u|)dx\right)
            \int_{\mathbb{R}^{N}}\phi_1(|\nabla u|)\nabla u\nabla \tilde{u}dx
 \nonumber\\
&&
       +    M_{2}\left(\int_{\mathbb{R}^{N}}\Phi_2(|\nabla v|)dx\right)
            \int_{\mathbb{R}^{N}}\phi_2(|\nabla v|)\nabla v\nabla \tilde{v}dx
       \nonumber\\
&&
          + \int_{\mathbb{R}^{N}}V_{1}(x)\phi_1(| u|) u \tilde{u}dx
          + \int_{\mathbb{R}^{N}}V_{2}(x)\phi_2(| v|) v \tilde{v}dx
\nonumber\\
&&        -\int_{\mathbb{R}^{N}}F_u(x,u,v)\tilde{u}dx
     -  \int_{\mathbb{R}^{N}}F_v(x,u,v)\tilde{v}dx
 \end{eqnarray}
 for all $(\tilde{u},\tilde{v})\in W$. Then the critical points of $I$ on $W$ are the weak solutions of system \eqref{eq1}.

\subsection{Existence}
\par
In this subsection, we present the following existence result and prove it by using the Mountain Pass Theorem.

\vskip2mm
\noindent
{\bf Theorem 3.1. }
  Assume that $(\phi_1)$--$(\phi_4)$, $(F_0)$, $(V_0)$, $(V_1)$, $(M_0)$, $(M_1)$ and the following conditions hold:
\begin{itemize}
	 \item[$(F_1)$] there exist two continuous functions $\psi_i ~(i=1, 2): [0,+\infty)\rightarrow \mathbb{R}$ and a constant $C_2>0$ such that
\begin{equation}\label{3.1.1}
 \left\{
  \begin{array}{l}
 |F_u(x,u,v)|\leq C_2\left(|u|^{l_{1}-1}+\psi_1(|u|)+\widetilde{\Psi}_1^{-1}(\Psi_2(|v|))\right),\\
 |F_v(x,u,v)|\leq C_2\left(|v|^{l_{2}-1}+\widetilde{\Psi}_2^{-1}(\Psi_1(|u|))+\psi_2(|v|)\right)
    \end{array}
 \right.
 \end{equation}
for all $(x,u,v)\in \mathbb{R}^{N} \times \mathbb{R}\times \mathbb{R}$,
where
$\Psi_i(t):=\int_{0}^{|t|}\psi_i(s)ds,$
 $t\in\mathbb{R} ~(i=1,2)$ are two $N$-functions satisfying
\begin{equation}\label{3.1.2}
       m_i
 <     l_{\Psi_{i}}:=\inf_{t>0}\frac{t\psi_i(t)}{\Psi_i(t)}
 \leq  \sup_{t>0}\frac{t\psi_i(t)}{\Psi_i(t)}
 =:    m_{\Psi_{i}}<l_{i}^{\ast},
 \end{equation}
 $\widetilde{\Psi}_i$ denote the complements of $\Psi_i ~(i=1,2)$, respectively;
 \item[$(F_2)$] there exists  a constant $C_3\in[0,1)$ such that
$$
\limsup_{|(u,v)|\rightarrow 0}\frac{F(x,u,v)}{\Phi_{1}(|u|)+\Phi_{2}(|v|)}=C_3 \mbox{\ \   uniformly in   } x \in \mathbb{R}^{N};
$$
\item[$(F_3)$] there exists a domain $G\subset \mathbb{R}^{N}$ such that
 $$
 \lim_{|(u,v)|\rightarrow+\infty}\frac{F(x,u,v)}{|u|^{m_1}+|v|^{m_2}}=+\infty, \mbox{\ \ for a.e.  } x \in G;
 $$
\item[$(F_4)$] there exist a continuous function $\overline{\gamma}: [0,\infty)\rightarrow \mathbb{R}^{+}$ and positive constants
$  \sigma_i\in
   \left[ \frac{l_i(m_{\overline{\Gamma}}-1)}{m_{\overline{\Gamma}}},
          \min\bigg\{l_i,
                    \frac{l_{i}^{\ast}(l_{\overline{\Gamma}}-1)}{l_{\overline{\Gamma}}}
              \bigg\}
   \right),$ $i=1,2,$
 $C_4, r>0$
 such that
 \begin{equation}\label{3.1.3}
      \overline{\Gamma}\left(\frac{F(x,u,v)}{|u|^{\sigma_1}+|v|^{\sigma_2}}\right)
 \leq C_4\overline{F}(x,u,v), \quad \mbox{ for all }x\in \mathbb{R}^{N} \mbox{ and }  (u,v)\in \mathbb R^2 \mbox{ with } |(u,v)|\geq r,
 \end{equation}
 where  $\overline{\Gamma}(t):=\int_{0}^{|t|}\overline{\gamma}(s)ds,~ t\in\mathbb{R}$, is an $N$-function with
 \begin{equation}\label{3.1.4}
       1
  <    l_{\overline{\Gamma}}:=\inf_{t>0}\frac{t\overline{\gamma}(t)}{\overline{\Gamma}(t)}
  \leq \sup_{t>0}\frac{t\overline{\gamma}(t)}{\overline{\Gamma}(t)}=:m_{\overline{\Gamma}}
  <+\infty
 \end{equation}
 and
 $$\overline{F}(x,u,v):=\frac{1}{m_1}F_u(x,u,v)u+\frac{1}{m_2}F_v(x,u,v)v-F(x,u,v)\geq 0, ~~~\forall  (x,u,v)\in \mathbb{R}^{N}\times \mathbb{R}\times \mathbb{R}.$$
\end{itemize}
 Then system \eqref{eq1} possesses a nontrivial weak solution.

\vskip2mm
 \noindent
{\bf Remark 3.2.}
 By  Lemma 2.7 and Lemma 3.1 in \cite{Bartsch2005},
the assumptions $(\phi_1)$--$(\phi_4)$, $(V_0)$, $(V_1)$ and (\ref{3.1.2}) imply that   the following embeddings are compact:
 $$
     W_{i}\hookrightarrow L^{\Psi_i}(\mathbb{R}^N),~
     W_{i}\hookrightarrow L^{\Phi_i}(\mathbb{R}^N)
      ~\mbox{ and }~
     W_{i}\hookrightarrow L^{p_i}(\mathbb{R}^N),
     ~~i=1,2,
 $$
  where $p_i \in [l_{i},l_{i}^{\ast})$.
 As a result, there exist some positive constants $C_{i,5},\;C_{i,6},\;i=1,2,$ such that
   \begin{eqnarray}\label{3.0.3}
       \|u\|_{L^{p_{i}}}\leq C_{i,5}\|u\|_{i},\;\|u\|_{L^{\Psi_i}}\leq C_{i,6}\|u\|_{i},
 \end{eqnarray}
  where $p_i \in [l_{i},l_{i}^{\ast})$.
 In particular, we have
 $\sigma_i\widetilde{l}_{\overline{\Gamma}},\;
 \sigma_i\widetilde{m}_{\overline{\Gamma}}\in [l_{i},l_{i}^{\ast})$, where $\sigma_i,i=1,2$,
$\widetilde{l}_{\overline{\Gamma}}=\frac{l_{\overline{\Gamma}}}{l_{\overline{\Gamma}}-1}$ and $\widetilde{m}_{\overline{\Gamma}}=\frac{m_{\overline{\Gamma}}}{m_{\overline{\Gamma}}-1}$. Hence, the following embeddings are compact:
 $$
     W_{i}\hookrightarrow L^{\sigma_i\widetilde{l}_{\overline{\Gamma}}}(\mathbb{R}^N)\;
     ~\mbox{ and }~
     W_{i}\hookrightarrow L^{\sigma_i\widetilde{m}_{\overline{\Gamma}}}(\mathbb{R}^N),
     ~~i=1,2.
 $$

 \noindent
{\bf Remark 3.3.}
 By  Young's inequality \eqref{2.1.1}, \eqref{3.1.1} and $F(x,0,0)=0$, the fact
 $$F(x,u,v)=\int_{0}^{u}F_s(x,s,v)ds+\int_{0}^{v}F_t(x,0,t)dt+F(x,0,0), \quad \forall  (x,u,v)\in\mathbb{R}^{N}\times\mathbb{R}\times\mathbb{R}$$
shows that there exists a constant $C_5>0$ such that
\begin{equation}\label{3.1.2-1}
|F(x,u,v)|\leq C_5(|u|^{l_{1}}+|v|^{l_{2}}+\Psi_1(|u|)+\Psi_2(|v|)), \quad \forall  (x,u,v)\in\mathbb{R}^{N}\times\mathbb{R}\times\mathbb{R}.
\end{equation}

 \noindent
{\bf Remark 3.4.} If we consider the system (\ref{eq1}) on a bounded domain $\Omega$ with Dirichlet boundary condition, then it is natural that we restrict those assumptions of Theorem 3.1 on the bounded domain $\Omega$. Thus we can claim that $(F_4)$ and $(f_2)$ are complementary.
 Firstly, we claim that if $l_{\Gamma}\ge m_{\overline{\Gamma}}$, $(F_4)$ and $(F_1)$ imply that the  condition $(f_2)$ hold. To be specific, choosing $r_{1}\geq \sqrt{2}r$ with $r>1$ and taking the suitable values of $l_{\Gamma}$ and $m_{\Gamma}$ such that
$1<  l_{\Gamma} \leq m_{\Gamma}<+\infty$,
$   \max\left\{1,\frac{(m_{\Psi_{1}}-l_1)m_{\Gamma}}{m_{\Psi_{1}}},\frac{(m_{\Psi_{2}}-l_2)m_{\Gamma}}{m_{\Psi_{2}}}\right\}
<    l_{\overline{\Gamma}}
\leq m_{\overline{\Gamma}}
<    l_{\Gamma}
\leq m_{\Gamma}$
and
 {\small $$
  m_{\Gamma}<    \min\left\{ \frac{l_{\overline{\Gamma}}m_{\overline{\Gamma}}(m_{\Psi_{1}}-l_1)+l_{\overline{\Gamma}}l_{1}}{m_{\overline{\Gamma}}(m_{\Psi_{1}}-l_1)},
                 \frac{l_1l_{\overline{\Gamma}}m_{\overline{\Gamma}}(m_{\Psi_{2}}-l_2)+l_{\overline{\Gamma}}l_{1}l_{2}}{l_1m_{\overline{\Gamma}}(m_{\Psi_{2}}-l_2)},
                 \frac{l_{\overline{\Gamma}}m_{\overline{\Gamma}}(m_{\Psi_{2}}-l_2)+l_{\overline{\Gamma}}l_{2}}{m_{\overline{\Gamma}}(m_{\Psi_{2}}-l_2)},
                 \frac{l_2l_{\overline{\Gamma}}m_{\overline{\Gamma}}(m_{\Psi_{1}}-l_1)+l_{\overline{\Gamma}}l_{1}l_{2}}{l_2m_{\overline{\Gamma}}(m_{\Psi_{1}}-l_1)}
         \right\},
$$}
 we can see that $(F_4)$ also hold  for $|(u,v)|\geq r_{1}$.
Obviously,
$(F_1)$ implies that
\begin{align}\label{3.1.2-11}
         |F(x,u,v)|
\leq   d_{2}(|u|^{l_{1}}+|v|^{l_{2}}+|u|^{m_{\Psi_{1}}}+|v|^{m_{\Psi_{2}}})\;\;\text { for }\;|(u,v)|\geq r,
 \end{align}
 where $d_{2}>0$.
Moreover, by \eqref{3.1.2-11} and Young's inequality, we have the following inequality
 \begin{align}\label{3.1.2-111}
&  (|u|^{\sigma_1}+|v|^{\sigma_2})^{l_{\overline{\Gamma}}}
       \left[d_{2}(|u|^{l_{1}}+|v|^{l_{2}}+|u|^{m_{\Psi_{1}}}+|v|^{m_{\Psi_{2}}})\right]^{m_{\Gamma}-l_{\overline{\Gamma}}}
      \nonumber\\
\leq &
         C_{m_{\Gamma}-l_{\overline{\Gamma}}}d_{2}^{m_{\Gamma}-l_{\overline{\Gamma}}}
         (|u|^{\sigma_1}+|v|^{\sigma_2})^{l_{\overline{\Gamma}}}
         (|u|^{l_{1}}+|v|^{l_{2}})^{m_{\Gamma}-l_{\overline{\Gamma}}}
         \nonumber\\
&~~~~
         +
      C_{l_{\overline{\Gamma}}}C^{2}_{m_{\Gamma}-l_{\overline{\Gamma}}}d_{2}^{m_{\Gamma}-l_{\overline{\Gamma}}}
         ( |u|^{\sigma_1l_{\overline{\Gamma}}+m_{\Psi_{1}}(m_{\Gamma}-l_{\overline{\Gamma}})}
          +|v|^{\sigma_2l_{\overline{\Gamma}}+m_{\Psi_{2}}(m_{\Gamma}-l_{\overline{\Gamma}})}
          )\nonumber\\
&~~~~
+
      C_{l_{\overline{\Gamma}}}C^{2}_{m_{\Gamma}-l_{\overline{\Gamma}}}d_{2}^{m_{\Gamma}-l_{\overline{\Gamma}}}
         \left(
            \frac{1}{\xi_{1}}|u|^{\sigma_1l_{\overline{\Gamma}}\xi_{1}}
          +\frac{\xi_{1}-1}{\xi_{1}}|v|^{\frac{\xi_{1}m_{\Psi_{2}}(m_{\Gamma}-l_{\overline{\Gamma}})}{\xi_{1}-1}}
          +\frac{\xi_{2}-1}{\xi_{2}}|u|^{\frac{\xi_{2}m_{\Psi_{1}}(m_{\Gamma}-l_{\overline{\Gamma}})}{\xi_{2}-1}}
          +\frac{1}{\xi_{2}}|v|^{\sigma_2l_{\overline{\Gamma}}\xi_{2}}
          \right)\nonumber\\
\leq &
         d_{3}
         \bigg[
         (|u|^{\sigma_1}+|v|^{\sigma_2})^{l_{\overline{\Gamma}}}
         (|u|^{l_{1}}+|v|^{l_{2}})^{m_{\Gamma}-l_{\overline{\Gamma}}}
        + |u|^{\sigma_1l_{\overline{\Gamma}}+m_{\Psi_{1}}(m_{\Gamma}-l_{\overline{\Gamma}})}
        + |v|^{\sigma_2l_{\overline{\Gamma}}+m_{\Psi_{2}}(m_{\Gamma}-l_{\overline{\Gamma}})}
        \nonumber\\
&~~~~
        + |u|^{\frac{\xi_{2}m_{\Psi_{1}}(m_{\Gamma}-l_{\overline{\Gamma}})}{\xi_{2}-1}}
        + |v|^{\frac{\xi_{1}m_{\Psi_{2}}(m_{\Gamma}-l_{\overline{\Gamma}})}{\xi_{1}-1}}
        + |u|^{\sigma_1l_{\overline{\Gamma}}\xi_{1}}
        + |v|^{\sigma_2l_{\overline{\Gamma}}\xi_{2}}
          \bigg],
\end{align}
for some
$\xi_{1}>\frac{l_2m_{\Gamma}}{l_2m_{\Gamma}-m_{\Psi_{2}}(m_{\Gamma}-l_{\overline{\Gamma}})}$,
$\xi_{2}>\frac{l_1m_{\Gamma}}{l_1m_{\Gamma}-m_{\Psi_{1}}(m_{\Gamma}-l_{\overline{\Gamma}})}$,
where
$
C_{m_{\Gamma}-l_{\overline{\Gamma}}}=
       \begin{cases}
              2^{m_{\Gamma}-l_{\overline{\Gamma}}-1},
                   \;\;\; \;\;\; \;\; \;             \text { if }\;\;m_{\Gamma}-l_{\overline{\Gamma}}>1, \\
              1,
                   \;\;\; \;\;\; \;\; \;             \text { if }\;\;m_{\Gamma}-l_{\overline{\Gamma}}\leq 1,
      \end{cases}
$
$C_{l_{\overline{\Gamma}}}=2^{l_{\overline{\Gamma}}-1}$,
$d_{3}>0$,
$
\sigma_1\in  \bigg[
           \frac{l_1(m_{\overline{\Gamma}}-1)}{m_{\overline{\Gamma}}},
           \min\bigg\{ l_1,
                             \frac{l_1m_{\Gamma}-m_{\Psi_{1}}(m_{\Gamma}-l_{\overline{\Gamma}})}{l_{\overline{\Gamma}}},
                             \frac{l_1l_2m_{\Gamma}-l_1m_{\Psi_{2}}(m_{\Gamma}-l_{\overline{\Gamma}})}{l_{\overline{\Gamma}}l_2},
                             \frac{l_{1}^{\ast}(l_{\overline{\Gamma}}-1)}{l_{\overline{\Gamma}}}
             \bigg\}
              \bigg)
$
and
$$
 \sigma_2 \in  \bigg[
          \frac{l_2(m_{\overline{\Gamma}}-1)}{m_{\overline{\Gamma}}},
          \min\bigg\{ l_2,
                             \frac{l_2m_{\Gamma}-m_{\Psi_{2}}(m_{\Gamma}-l_{\overline{\Gamma}})}{l_{\overline{\Gamma}}},
                             \frac{l_2l_1m_{\Gamma}-l_2m_{\Psi_{1}}(m_{\Gamma}-l_{\overline{\Gamma}})}{l_{\overline{\Gamma}}l_1},
                             \frac{l_{2}^{\ast}(l_{\overline{\Gamma}}-1)}{l_{\overline{\Gamma}}}
               \bigg\}
                 \bigg)
$$
 which is  reasonable by
our example in the last section with $m_{\Gamma}=2$.
Hence,
in virtue of  $(F_4)$, \eqref{3.1.2-11}, \eqref{3.1.2-111}, Young's inequality  and Lemma 2.3, we have
 \begin{align}\label{3.1.3-}
&     C_4\overline{F}(x,u,v)
\nonumber\\
&\geq \overline{\Gamma}\left(\frac{F(x,u,v)}{|u|^{\sigma_1}+|v|^{\sigma_2}}\right)
  \nonumber\\
&\geq \overline{\Gamma}(1)
\frac{
     \min\left\{
                 \left(\frac{|F(x,u,v)|}{|u|^{\sigma_1}+|v|^{\sigma_2}}\right)^{l_{\overline{\Gamma}}},
                 \left(\frac{|F(x,u,v)|}{|u|^{\sigma_1}+|v|^{\sigma_2}}\right)^{m_{\overline{\Gamma}}}
         \right\}
     }
     {
     \max\left\{
                 \left(\frac{|F(x,u,v)|}{|u|^{l_1}+|v|^{l_2}}\right)^{l_{\Gamma}},
                 \left(\frac{|F(x,u,v)|}{|u|^{l_1}+|v|^{l_2}}\right)^{m_{\Gamma}}
         \right\}
     }
     \max\bigg\{
                 \left(\frac{|F(x,u,v)|}{|u|^{l_1}+|v|^{l_2}}\right)^{l_{\Gamma}},
                 \left(\frac{|F(x,u,v)|}{|u|^{l_1}+|v|^{l_2}}\right)^{m_{\Gamma}}
         \bigg\}
         \nonumber\\
&\geq \begin{cases}
          \frac{\overline{\Gamma}(1)}{\Gamma(1)}
     \Gamma\left(\frac{|F(x,u,v)|}{|u|^{l_1}+|v|^{l_2}}\right)
        \frac{(|u|^{l_1}+|v|^{l_2})^{m_{\Gamma}}
             }
             {
               (|u|^{\sigma_1}+|v|^{\sigma_2})^{m_{\overline{\Gamma}}}
               |F(x,u,v)|^{m_{\Gamma}-m_{\overline{\Gamma}}}
             }, \;\;\; \;\;\; \;\; \;             \text { if }\;\;\frac{|F(x,u,v)|}{|u|^{l_1}+|v|^{l_2}}\geq 1,
                                                            \;\frac{|F(x,u,v)|}{|u|^{\sigma_1}+|v|^{\sigma_2}}< 1, \\
          \frac{\overline{\Gamma}(1)}{\Gamma(1)}
     \Gamma\left(\frac{|F(x,u,v)|}{|u|^{l_1}+|v|^{l_2}}\right)
        \frac{(|u|^{l_1}+|v|^{l_2})^{m_{\Gamma}}
             }
             {
               (|u|^{\sigma_1}+|v|^{\sigma_2})^{l_{\overline{\Gamma}}}
               |F(x,u,v)|^{m_{\Gamma}-l_{\overline{\Gamma}}}
             }, \;\;\; \;\;\; \;\; \;             \text { if }\;\;\frac{|F(x,u,v)|}{|u|^{l_1}+|v|^{l_2}}\geq 1,
                                                             \;\frac{|F(x,u,v)|}{|u|^{\sigma_1}+|v|^{\sigma_2}}\geq 1, \\
     \frac{\overline{\Gamma}(1)}{\Gamma(1)}
     \Gamma\left(\frac{|F(x,u,v)|}{|u|^{l_1}+|v|^{l_2}}\right)
        \frac{(|u|^{l_1}+|v|^{l_2})^{l_{\Gamma}}
             }
             {
             (|u|^{\sigma_1}+|v|^{\sigma_2})^{m_{\overline{\Gamma}}}
             |F(x,u,v)|^{l_{\Gamma}-m_{\overline{\Gamma}}}
             }, \;\;\; \;\;\; \;\; \;             \text { if }\;\;\frac{|F(x,u,v)|}{|u|^{l_1}+|v|^{l_2}}< 1,
                                                        \;\frac{|F(x,u,v)|}{|u|^{\sigma_1}+|v|^{\sigma_2}}< 1,\\
     \frac{\overline{\Gamma}(1)}{\Gamma(1)}
     \Gamma\left(\frac{|F(x,u,v)|}{|u|^{l_1}+|v|^{l_2}}\right)
        \frac{(|u|^{l_1}+|v|^{l_2})^{l_{\Gamma}}
        }
        {
        (|u|^{\sigma_1}+|v|^{\sigma_2})^{l_{\overline{\Gamma}}}
         |F(x,u,v)|^{l_{\Gamma}-l_{\overline{\Gamma}}}
        }, \;\;\; \;\;\; \;\; \;             \text { if }\;\;\frac{|F(x,u,v)|}{|u|^{l_1}+|v|^{l_2}}< 1,
                                                        \;\frac{|F(x,u,v)|}{|u|^{\sigma_1}+|v|^{\sigma_2}}\geq 1
          \end{cases}
          \nonumber\\
&\geq \begin{cases}
          \frac{\overline{\Gamma}(1)}{\Gamma(1)}
     \Gamma\left(\frac{|F(x,u,v)|}{|u|^{l_1}+|v|^{l_2}}\right)
        \frac{(|u|^{l_1}+|v|^{l_2})^{m_{\Gamma}}
             }
             {
               (|u|^{\sigma_1}+|v|^{\sigma_2})^{m_{\overline{\Gamma}}}
               (|u|^{\sigma_1}+|v|^{\sigma_2})^{m_{\Gamma}-m_{\overline{\Gamma}}}
             }, \;\;\; \;\;\; \;\; \;             \text { if }\;\;\frac{|F(x,u,v)|}{|u|^{l_1}+|v|^{l_2}}\geq 1,
                                                            \;\frac{|F(x,u,v)|}{|u|^{\sigma_1}+|v|^{\sigma_2}}< 1, \\
          \frac{\overline{\Gamma}(1)}{\Gamma(1)}
     \Gamma\left(\frac{|F(x,u,v)|}{|u|^{l_1}+|v|^{l_2}}\right)
        \frac{(|u|^{l_1}+|v|^{l_2})^{m_{\Gamma}}
             }
             {
               (|u|^{\sigma_1}+|v|^{\sigma_2})^{l_{\overline{\Gamma}}}
               \left[d_{2}(|u|^{l_{1}}+|v|^{l_{2}}+|u|^{m_{\Psi_{1}}}+|v|^{m_{\Psi_{2}}})
               \right]^{m_{\Gamma}-l_{\overline{\Gamma}}}
             }, \;\;             \text { if }\;\;\frac{|F(x,u,v)|}{|u|^{l_1}+|v|^{l_2}}\geq 1,
                                                             \;\frac{|F(x,u,v)|}{|u|^{\sigma_1}+|v|^{\sigma_2}}\geq 1, \\
     \frac{\overline{\Gamma}(1)}{\Gamma(1)}
     \Gamma\left(\frac{|F(x,u,v)|}{|u|^{l_1}+|v|^{l_2}}\right)
        \frac{(|u|^{l_1}+|v|^{l_2})^{l_{\Gamma}}
             }
             {
             (|u|^{\sigma_1}+|v|^{\sigma_2})^{m_{\overline{\Gamma}}}
             (|u|^{l_1}+|v|^{l_2})^{l_{\Gamma}-m_{\overline{\Gamma}}}
             }, \;\;\; \;\;\; \;\; \;             \text { if }\;\;\frac{|F(x,u,v)|}{|u|^{l_1}+|v|^{l_2}}< 1,
                                                        \;\frac{|F(x,u,v)|}{|u|^{\sigma_1}+|v|^{\sigma_2}}< 1,\\
    \frac{\overline{\Gamma}(1)}{\Gamma(1)}
     \Gamma\left(\frac{|F(x,u,v)|}{|u|^{l_1}+|v|^{l_2}}\right)
        \frac{(|u|^{l_1}+|v|^{l_2})^{l_{\Gamma}}
        }
        {
        (|u|^{\sigma_1}+|v|^{\sigma_2})^{l_{\overline{\Gamma}}}
         (|u|^{l_1}+|v|^{l_2})^{l_{\Gamma}-l_{\overline{\Gamma}}}
        }, \;\;\; \;\;\; \;\; \;             \text { if }\;\;\frac{|F(x,u,v)|}{|u|^{l_1}+|v|^{l_2}}< 1,
                                                        \;\frac{|F(x,u,v)|}{|u|^{\sigma_1}+|v|^{\sigma_2}}\geq 1
          \end{cases}
          \nonumber\\
&\geq d_{4}\frac{\overline{\Gamma}(1)}{\Gamma(1)}
     \Gamma\left(\frac{|F(x,u,v)|}{|u|^{l_1}+|v|^{l_2}}\right)
 \end{align}
for $|(u,v)|\geq r_{1}$,  where
$d_{4}>0$.
To be specific, for $|u|, |v|\geq r$, we have
\begin{align}\label{3.1.3-1}
   \frac{(|u|^{l_1}+|v|^{l_2})^{m_{\Gamma}}
        }
        {(|u|^{\sigma_1}+|v|^{\sigma_2})^{m_{\overline{\Gamma}}}
         (|u|^{\sigma_1}+|v|^{\sigma_2})^{m_{\Gamma}-m_{\overline{\Gamma}}}
        }
\geq
     \min\left\{r^{m_{\Gamma}(l_1-\sigma_1)},r^{m_{\Gamma}(l_2-\sigma_2)}\right\},
\end{align}
\begin{align}\label{3.1.3-2}
    \frac{(|u|^{l_1}+|v|^{l_2})^{l_{\Gamma}}
             }
             {
             (|u|^{\sigma_1}+|v|^{\sigma_2})^{m_{\overline{\Gamma}}}
             (|u|^{l_1}+|v|^{l_2})^{l_{\Gamma}-m_{\overline{\Gamma}}}
             }
\geq
     \min\left\{r^{m_{\overline{\Gamma}}(l_1-\sigma_1)},r^{m_{\overline{\Gamma}}(l_2-\sigma_2)}\right\},
\end{align}
\begin{align}\label{3.1.3-3}
&   \frac{(|u|^{l_1}+|v|^{l_2})^{l_{\Gamma}}
        }
        {
        (|u|^{\sigma_1}+|v|^{\sigma_2})^{l_{\overline{\Gamma}}}
         (|u|^{l_1}+|v|^{l_2})^{l_{\Gamma}-l_{\overline{\Gamma}}}
        }
\geq
     \min\left\{r^{l_{\overline{\Gamma}}(l_1-\sigma_1)},r^{l_{\overline{\Gamma}}(l_2-\sigma_2)}\right\},
\end{align}
\begin{align}\label{3.1.3-4}
    \frac{(|u|^{\sigma_1}+|v|^{\sigma_2})^{l_{\overline{\Gamma}}}
          (|u|^{l_{1}}+|v|^{l_{2}})^{m_{\Gamma}-l_{\overline{\Gamma}}}
         }
         {
          (|u|^{l_1}+|v|^{l_2})^{m_{\Gamma}}
         }
\leq
     \frac{1}{\max\left\{r^{l_{\overline{\Gamma}}(l_1-\sigma_1)},r^{l_{\overline{\Gamma}}(l_2-\sigma_2)}\right\}},
\end{align}
\begin{align}\label{3.1.3-5}
    \frac{ |u|^{\sigma_1l_{\overline{\Gamma}}+m_{\Psi_{1}}(m_{\Gamma}-l_{\overline{\Gamma}})}
          +|v|^{\sigma_2l_{\overline{\Gamma}}+m_{\Psi_{2}}(m_{\Gamma}-l_{\overline{\Gamma}})}
         }
         {
          (|u|^{l_1}+|v|^{l_2})^{m_{\Gamma}}
         }
\leq \frac{1}{\max\left\{r^{l_1m_{\Gamma}-\sigma_1l_{\overline{\Gamma}}-m_{\Psi_{1}}(m_{\Gamma}-l_{\overline{\Gamma}})},
                          r^{l_2m_{\Gamma}-\sigma_2l_{\overline{\Gamma}}-m_{\Psi_{2}}(m_{\Gamma}-l_{\overline{\Gamma}})}
                   \right\}
              },
\end{align}
\begin{align}\label{3.1.3-6}
    \frac{ |u|^{\frac{\xi_{2}m_{\Psi_{1}}(m_{\Gamma}-l_{\overline{\Gamma}})}{\xi_{2}-1}}
          +|v|^{\frac{\xi_{1}m_{\Psi_{2}}(m_{\Gamma}-l_{\overline{\Gamma}})}{\xi_{1}-1}}
         }
         {
          (|u|^{l_1}+|v|^{l_2})^{m_{\Gamma}}
         }
\leq
      \frac{1}{
                  \max\left\{r^{l_1m_{\Gamma}-\frac{\xi_{2}m_{\Psi_{1}}(m_{\Gamma}-l_{\overline{\Gamma}})}{\xi_{2}-1}},
                              r^{l_2m_{\Gamma}-\frac{\xi_{1}m_{\Psi_{2}}(m_{\Gamma}-l_{\overline{\Gamma}})}{\xi_{1}-1}}
                     \right\}
              }
\end{align}
and
\begin{align}\label{3.1.3-7}
    \frac{ |u|^{\sigma_1l_{\overline{\Gamma}}\xi_{1}}
          +|v|^{\sigma_2l_{\overline{\Gamma}}\xi_{2}}
         }
         {
          (|u|^{l_1}+|v|^{l_2})^{m_{\Gamma}}
         }
\leq
      \frac{1}{
                  \max\left\{
                              r^{l_1m_{\Gamma}-\sigma_1l_{\overline{\Gamma}}\xi_{1}},
                              r^{l_2m_{\Gamma}-\sigma_2l_{\overline{\Gamma}}\xi_{2}}
                     \right\}
              }.
\end{align}
For $|u|\geq r, 0 \leq |v|<r$ and $|v|\geq r, 0 \leq |u|<r$, the inequalities \eqref{3.1.3-1}-\eqref{3.1.3-7} also hold with different values on the right side of the inequalities.
Hence, \eqref{3.1.3-} holds for $|(u,v)|\geq r_{1}$ and then it is easy to  see that $(f_{2})$ holds.
\par
Next, we claim that $(f_2)$ and the following assumption (A) imply that  condition $(F_4)$ holds if we take
 $m_{\Gamma}\geq l_{\Gamma}\geq \frac{m_{i}m_{\overline{\Gamma}}}{m_{i}-l_{i}m_{\overline{\Gamma}}}$.
\par{\it (A)
 $\overline{F}(x,u,v)\geq C (|u|^{\varrho_{1}}+|v|^{\varrho_{2}})$ for all $|(u,v)|\geq r_{1}$, where
 $\varrho_i \in (0, m_{i}]$ with $m_{i}> l_{i}m_{\overline{\Gamma}}$.
 }
 \par
 In fact,
 choosing $r= \sqrt{2}r_{1}$ with $r_{1} >1$,
we can see that $(f_2)$ also holds  for $(u,v)\in \mathbb{R}^2$ with $|(u,v)|\geq r$ and $|u|^{\sigma_1}+|v|^{\sigma_2}\geq 1$. Next,  there are four possible cases, that is,
\textcircled{1}
$\frac{|F(x,u,v)|}{|u|^{l_1}+|v|^{l_2}}<1, \frac{|F(x,u,v)|}{|u|^{\sigma_1}+|v|^{\sigma_2}}\geq 1$,
$\varrho_{i}\geq \frac{l_{i}m_{\Gamma}m_{\overline{\Gamma}}}{m_{\Gamma}-m_{\overline{\Gamma}}}$;
\textcircled{2}
$\frac{|F(x,u,v)|}{|u|^{l_1}+|v|^{l_2}}<1, \frac{|F(x,u,v)|}{|u|^{\sigma_1}+|v|^{\sigma_2}}< 1$,
$\varrho_{i}\geq \frac{l_{i}m_{\Gamma}l_{\overline{\Gamma}}}{m_{\Gamma}-l_{\overline{\Gamma}}}$;
\textcircled{3}
$\frac{|F(x,u,v)|}{|u|^{l_1}+|v|^{l_2}} \geq 1$, $\frac{|F(x,u,v)|}{|u|^{\sigma_1}+|v|^{\sigma_2}}\geq 1$,
$\varrho_{i}\geq \frac{l_{i}l_{\Gamma}m_{\overline{\Gamma}}}{l_{\Gamma}-m_{\overline{\Gamma}}}$
and
\textcircled{4}
$\frac{|F(x,u,v)|}{|u|^{l_1}+|v|^{l_2}} \geq 1$, $\frac{|F(x,u,v)|}{|u|^{\sigma_1}+|v|^{\sigma_2}} < 1$,
$\varrho_{i}\geq \frac{l_{i}l_{\Gamma}l_{\overline{\Gamma}}}{l_{\Gamma}-l_{\overline{\Gamma}}}$.
Without loss of generality, we only focus on the proof of the first case.
\par
\noindent
\textcircled{1} Since $\frac{|F(x,u,v)|}{|u|^{l_1}+|v|^{l_2}}<1, \frac{|F(x,u,v)|}{|u|^{\sigma_1}+|v|^{\sigma_2}}\geq 1$
and $\varrho_{i}\geq \frac{l_{i}m_{\Gamma}m_{\overline{\Gamma}}}{m_{\Gamma}-m_{\overline{\Gamma}}}$, by $(f_2)$, (A) and Lemma 2.3,
we have
\begin{align*}
 d_{1} \overline{F}(x,u,v)
&= d_{1} \left(\overline{F}(x,u,v)\right)^{\frac{m_{\overline{\Gamma}}}{m_{\Gamma}}}
             \left(\overline{F}(x,u,v)\right)^{\frac{m_{\Gamma}-m_{\overline{\Gamma}}}{m_{\Gamma}}}\\
&\geq    d_{1}^{\frac{m_{\Gamma}-m_{\overline{\Gamma}}}{m_{\Gamma}}}
         \left(\Gamma\left(\frac{F(x,u,v)}{|u|^{l_1}+|v|^{l_2}}\right)\right)^{\frac{m_{\overline{\Gamma}}}{m_{\Gamma}}}
         \left(\overline{F}(x,u,v)\right)^{\frac{m_{\Gamma}-m_{\overline{\Gamma}}}{m_{\Gamma}}}
  \\
&\geq    d_{1}^{\frac{m_{\Gamma}-m_{\overline{\Gamma}}}{m_{\Gamma}}}
         \Gamma(1)
         \left(\frac{|F(x,u,v)|}{|u|^{l_1}+|v|^{l_2}}\right)^{m_{\overline{\Gamma}}}
         \left(|u|^{\varrho_1}+|v|^{\varrho_2}\right)^{\frac{m_{\Gamma}-m_{\overline{\Gamma}}}{m_{\Gamma}}}
  \\
& =     d_{1}^{\frac{m_{\Gamma}-m_{\overline{\Gamma}}}{m_{\Gamma}}}
        \Gamma(1)
       \left(\frac{|F(x,u,v)|}{|u|^{\sigma_1}+|v|^{\sigma_2}}\right)^{m_{\overline{\Gamma}}}
        \left(|u|^{\sigma_1}+|v|^{\sigma_2}\right)^{m_{\overline{\Gamma}}}
        \frac{
              \left(|u|^{\varrho_1}+|v|^{\varrho_2}\right)^{\frac{m_{\Gamma}-m_{\overline{\Gamma}}}{m_{\Gamma}}}
             }
             {\left(|u|^{l_1}+|v|^{l_2}\right)^{m_{\overline{\Gamma}}}
             }
  \\
& \geq   d_{1}^{\frac{m_{\Gamma}-m_{\overline{\Gamma}}}{m_{\Gamma}}}
           \frac{\Gamma(1)}{\overline{\Gamma}(1)}
          \overline{\Gamma}\left(\frac{|F(x,u,v)|}{|u|^{\sigma_1}+|v|^{\sigma_2}}\right)
          \frac{
              \left(|u|^{\varrho_1}+|v|^{\varrho_2}\right)^{\frac{m_{\Gamma}-m_{\overline{\Gamma}}}{m_{\Gamma}}}
             }
             {\left(|u|^{l_1}+|v|^{l_2}\right)^{m_{\overline{\Gamma}}}
             }
  \\
& \geq
          \overline{\Gamma}\left(\frac{|F(x,u,v)|}{|u|^{\sigma_1}+|v|^{\sigma_2}}\right)
          \frac{d_{1}^{\frac{m_{\Gamma}-m_{\overline{\Gamma}}}{m_{\Gamma}}}\Gamma(1)}{C_{m_{\overline{\Gamma}}}\overline{\Gamma}(1)}
          \left( \frac{    |u|^{\frac{\varrho_1(m_{\Gamma}-m_{\overline{\Gamma}})}{m_{\Gamma}}}
                         + |v|^{\frac{\varrho_2(m_{\Gamma}-m_{\overline{\Gamma}})}{m_{\Gamma}}}
                      }
                      {    |u|^{l_1m_{\overline{\Gamma}}}
                         + |v|^{l_2m_{\overline{\Gamma}}}
                      }
          \right)
  \\
& \geq    d_{5}
          \overline{\Gamma}\left(\frac{|F(x,u,v)|}{|u|^{\sigma_1}+|v|^{\sigma_2}}\right),
\end{align*}
where $d_{5}>0$ and $C_{m_{\overline{\Gamma}}}=2^{m_{\overline{\Gamma}}-1}$.
\par
Hence, $(F_4)$ and $(f_2)$ are complementary. We will give an example which satisfies $(F_4)$ but not satisfies $(f_2)$  in Section 5.

 \vskip2mm
 \noindent
{\bf Lemma 3.5.}
  Suppose that $(\phi_1)$--$(\phi_3)$, $(M_0)$, $(V_0)$ and $(F_1)$  hold. Then there are constants $\rho, \alpha$$ >0$ such that $I\mid _{\partial B_{\rho}}\geq \alpha$.

\vskip2mm
\noindent
{\bf Proof.}
 By $(\phi_1)$--$(\phi_3)$ and $(V_0)$,  we have
 \begin{equation}\label{3.1.6}
        \min\left\{\|u\|_{\Phi_1,V_1}^{l_1},\|u\|_{\Phi_1,V_1}^{m_1}\right\}
 \leq   \int_{\mathbb{R}^{N}}V_{1}(x)\Phi_1(| u|)dx
 \leq   \max\left\{\|u\|_{\Phi_1,V_1}^{l_1},\|u\|_{\Phi_1,V_1}^{m_1}\right\}
 \end{equation}
 and
 \begin{equation}\label{3.1.7}
 \min\left\{\|v\|_{\Phi_2,V_2}^{l_2},\|v\|_{\Phi_2,V_2}^{m_2}\right\}
 \leq   \int_{\mathbb{R}^N}V_{2}(x)\Phi_2(| v|)dx
 \leq   \max\left\{\|v\|_{\Phi_2,V_2}^{l_2},\|v\|_{\Phi_2,V_2}^{m_2}\right\}
 \end{equation}
 for all $u\in W_{1}$ and  $v\in W_{2}$ (the details is in Lemma 2.1 of \cite{Liu-Shibo2019}).
 Moreover, by \eqref{3.1.2-1} and $(F_{2})$,
 there exist  constants $C_{6}\in (0,1)$ and $C_{7}>0$ such that
 \begin{equation}\label{3.1.5+}
 |F(x,u,v)|\leq  (1-C_{6})(\Phi_1(|u|)+\Phi_2(|v|))+C_{7}(\Psi_1(|u|)+\Psi_2(|v|)), \quad \forall  (x,u,v)\in \mathbb{R}^N\times \mathbb{R}\times \mathbb{R}.
 \end{equation}
Choosing $\rho>0$ such that $\rho=\|(u,v)\|=\|u\|_{1}+\|v\|_{2} <\min\left\{\frac{1}{\max\{C_{1,6},C_{2,6}\}},1\right\}$. Then
$\|u\|_{\Psi_1}\leq C_{1,6}\|u\|_{1}<1$
and $\|v\|_{\Psi_2}\leq C_{2,6}\|v\|_{2}<1$.
 By  \eqref{3.1.6}, \eqref{3.1.7}, \eqref{3.1.5+}, $(M_{0})$, Lemma 2.4
and Remark 3.2, we obtain
 \begin{align*}
      I(u,v)
&\geq       C_{1,3}\int_{\mathbb{R}^N}\Phi_1(|\nabla u|)dx
          + \int_{\mathbb{R}^N}V_{1}(x)\Phi_1(| u|)dx
          - (1-C_{6})\int_{\mathbb{R}^N}V_{1}(x)\Phi_{1}(u)dx
          - C_{7}\int_{\mathbb{R}^N}\Psi_1(u)dx
          \nonumber\\
&~~~~
          + C_{2,3}\int_{\mathbb{R}^N}\Phi_2(|\nabla v|)dx
          + \int_{\mathbb{R}^N}V_{2}(x)\Phi_2(| v|)dx
          - (1-C_{6})\int_{\mathbb{R}^N}V_{2}(x)\Phi_{2}(v)dx
          - C_{7}\int_{\mathbb{R}^N}\Psi_2(v)dx
           \\
&\geq    C_{1,3}\min\left\{\|\nabla u\|_{\Phi_1}^{l_1},\|\nabla u\|_{\Phi_1}^{m_1}\right\}
       + C_{6}\min\left\{\|u\|_{\Phi_1,V_1}^{l_1},\|u\|_{\Phi_1,V_1}^{m_1}\right\}
       - C_{7}\max\left\{\|u\|_{\Psi_1}^{l_{\Psi_1}},\|u\|_{\Psi_1}^{m_{\Psi_1}}\right\}
       \\
&~~~
       + C_{2,3}\min\left\{\|\nabla v\|_{\Phi_2}^{l_2},\|\nabla v\|_{\Phi_2}^{m_2}\right\}
       + C_{6}\min\left\{\|v\|_{\Phi_2,V_2}^{l_2},\|v\|_{\Phi_2,V_2}^{m_2}\right\}
       - C_{7}\max\left\{\|v\|_{\Psi_2}^{l_{\Psi_2}},\|v\|_{\Psi_2}^{m_{\Psi_2}}\right\}
       \\
&\geq    C_{1,3}\|\nabla u\|_{\Phi_1}^{m_1}
       + C_{6}\|u\|_{\Phi_1,V_1}^{m_1}
       - C_{7}C_{1,6}\|u\|_{1}^{l_{\Psi_1}}
       + C_{2,3}\|\nabla v\|_{\Phi_2}^{m_2}
       + C_{6}\|v\|_{\Phi_2,V_2}^{m_2}
       - C_{7}C_{2,6}\|v\|_{2}^{l_{\Psi_2}}
       \\
&\geq   \|u\|_{1}^{m_1}\left(  \frac{\min\left\{C_{1,3},C_{6}\right\}}{2^{m_1-1}}
                                      - C_{7}C_{1,6}\|u\|_{1}^{l_{\Psi_1}-m_1}
                               \right)
       + \|v\|_{2}^{m_2}\left(   \frac{\min\left\{C_{2,3},C_{6}\right\}}{2^{m_2-1}}
                                      - C_{7}C_{2,6}\|v\|_{2}^{l_{\Psi_2}-m_2}
                               \right).
\end{align*}
Since $1<m_i<l_{\Psi_{i}}$, we can choose positive constants $\rho$ and $\alpha$ small enough such that ${I}(u,v)\geq \alpha$ for all $(u,v)\in W$ with $\|(u,v)\|=\rho$. \qed

\vskip2mm
\noindent
 {\bf Lemma 3.6.}
 Suppose that $(\phi_1)$--$(\phi_3)$,  $(M_0)$, $(V_0)$ and $(F_3)$ hold. Then there is a point $(u,v)\in W\backslash B_{\rho}$ such that $I(u,v)\leq0$.
 \vskip2mm
 \noindent
{\bf Proof.}   By $(F_3)$ and the continuity of $F$, there exist two constants $C_{8}>0$ and $C_{9}>0$ such that
\begin{equation}\label{3.1.8}
 F(x,u,v)\geq C_{8}(|u|^{m_1}+|v|^{m_2})-C_{9},\;\forall  (u,v)\in \mathbb{R}\times \mathbb{R}\; \mbox{and\;a.e.}\;x\in G.
 \end{equation}
Choose $u_0\in C_{0}^{\infty}(\mathbb{R}^N)\setminus\{0\}$ with $0< u_0(x)\leq 1$ and $supp(u_0) \subset\subset G$. Obviously, $(tu_0,0)\in W$ for all $t\in\mathbb{R}$.   By  $(M_0)$,  (2) in Lemma 2.3, \eqref{3.1.6}, \eqref{3.1.7} and \eqref{3.1.8},   when $t>1$, we have
 \begin{align*}
            I(tu_0,0)
 &   =      \widehat{M_{1}}\left(\int_{\mathbb{R}^N}\Phi_1(t|\nabla u_0|)dx\right)
          + \int_{\mathbb{R}^N}V_{1}(x)\Phi_1(t| u_0|)dx
          - \int_{G}F(x,tu_0,0)dx\\
 &\leq
            C_{1,4}t^{m_1}\int_{\mathbb{R}^N}\Phi_1(|\nabla u_0|)dx
          + t^{m_1}\int_{\mathbb{R}^N}V_{1}(x)\Phi_1(| u_0|)dx
          -  C_{8}\int_{G}t^{m_1}|u_0|^{m_1}dx
          +  C_{9} |\mbox{supp}u_0| \\
 &\leq     t^{m_1}
            \bigg(C_{1,4}\|\nabla u_0\|_{\Phi_1}^{l_1}
                 +C_{1,4}\|\nabla u_0\|_{\Phi_1}^{m_1}
                 +\|u_0\|_{\Phi_{1},V_{1}}^{l_1}
                 +\|u_0\|_{\Phi_{1},V_{1}}^{m_1}
                 -C_{8}\|u_0\|_{L^{m_1}(G)}^{m_1}
             \bigg)
             +  C_{9} |\mbox{supp}u_0|.
 \end{align*}
If we choose
\begin{align*}
  C_{8}>
  \frac{   C_{1,4}\||\nabla u_0|\|_{\Phi_1}^{l_1}
         + C_{1,4}\||\nabla u_0|\|_{\Phi_1}^{m_1}
         + \|u_0\|_{\Phi_{1},V_{1}}^{l_1}
         + \|u_0\|_{\Phi_{1},V_{1}}^{m_1}
      }{\|u_0\|_{L^{m_1}(G)}^{m_1}},
\end{align*}
then there exists $t$ large enough such that
 $I(tu_0,0)\leq 0$ and $\|(tu_0,0)\|>\rho.$ \qed

\vskip2mm
\noindent
  {\bf Lemma 3.7.}
   Suppose that $(\phi_1)$--$(\phi_4)$, $(V_0)$,  $(V_1)$,  $(M_1)$,  $(F_1)$, $(F_3)$ and $(F_4)$ hold. Then (C)-sequence in $W$  is bounded.
   \vskip2mm
  \noindent
{\bf Proof.}  Let $\{(u_n,v_n)\}$ be a $(C)$-sequence of $I$ in $W$. Then, for $n$ large enough, by $(\phi_3)$,  $(M_1)$ and \eqref{2.2.1}, there exists a $c>0$ such that
 \begin{align}\label{3.1.10}
          c+1
 &\geq    I(u_n,v_n)
         -\left\langle I'(u_n,v_n),\left(\frac{1}{m_1}u_n,\frac{1}{m_2}v_n\right)\right\rangle\nonumber\\
 & =        \widehat{M_{1}}\left(\int_{\mathbb{R}^N}\Phi_1(|\nabla u_n|)dx\right)
        - \frac{1}{m_1}
            M_{1}\left(\int_{\mathbb{R}^N}\Phi_1(|\nabla u_n|)dx\right)
                   \int_{\mathbb{R}^N}\phi_1(|\nabla u_n|)|\nabla u_n|^{2}dx\nonumber\\
&~~~
        + \widehat{M_{2}}\left(\int_{\mathbb{R}^N}\Phi_2(|\nabla v_n|)dx\right)
        -  \frac{1}{m_2}
            M_{2}\left(\int_{\mathbb{R}^N}\Phi_2(|\nabla v_n|)dx\right)
                   \int_{\mathbb{R}^N}\phi_2(|\nabla v_n|)|\nabla v_n|^{2}dx\nonumber\\
&~~~
          + \int_{\mathbb{R}^N}\left(V_{1}(x)\Phi_1(|u_n|)
                                     - \frac{1}{m_1}V_{1}(x)\phi_1(|u_n|) |u_n|^{2}
                                  \right)dx
          \nonumber\\
&~~~
          + \int_{\mathbb{R}^N}\left(V_{2}(x)\Phi_2(|v_n|)
                                      - \frac{1}{m_2}V_{2}(x)\phi_2(|v_n|) |v_n|^{2}
                                    \right)dx
\nonumber\\
&~~~   +\int_{\mathbb{R}^N}\left(\frac{1}{m_1}F_u(x,u_n,v_n)u_n
                               +\frac{1}{m_2}F_v(x,u_n,v_n)v_n
                               -F(x,u_n,v_n)
                               \right)dx\nonumber\\
&\geq  \int_{\mathbb{R}^N}\overline{F}(x,u_n,v_n)dx.
\end{align}
 Arguing by contradiction, we assume that there exists a subsequence of $\{(u_n,v_n)\}$, still denoted by $\{(u_n,v_n)\}$, such that $\|(u_n,v_n)\|=\|u_n\|_{1}+\|v_n\|_{2}\rightarrow +\infty$. Then we discuss this problem  in two situations.
 \vskip2mm
 \noindent {\bf Case 1}. Suppose that $\|u_n\|_{1}\rightarrow +\infty$ and
                           $\|v_n\|_{2}\rightarrow +\infty$.
 Let $\bar{u}_n=\frac{u_n}{\|u_n\|_{1}}$ and
     $\bar{v}_n=\frac{v_n}{\|v_n\|_{2}}$.
Then $\{(\bar{u}_n,\bar{v}_n)\}$ is bounded in $W$.
Passing to a subsequence $\{(\bar{u}_n,\bar{v}_n)\}$, by Remark 3.2, there exists a point $(\bar{u},\bar{v})\in W$ such that\\
 $\star$\quad $\bar{u}_n\rightharpoonup \bar{u}$ in $W_{1}$,
\; $\bar{u}_n\rightarrow \bar{u}$
                                          in $L^{l_1}(\mathbb{R}^{N})$,
                                          in $L^{\sigma_1\widetilde{l_{\overline{\Gamma}}}}(\mathbb{R}^{N})$
                                      and in $L^{\sigma_1\widetilde{m_{\overline{\Gamma}}}}(\mathbb{R}^{N})$,
 $\bar{u}_n(x)\rightarrow \bar{u}(x)$ a.e. in $\mathbb{R}^{N}$;\\
 $\star$\quad $\bar{v}_n\rightharpoonup \bar{v}$ in $W_{2}$,
 \quad $\bar{v}_n\rightarrow \bar{v}$
                                            in $L^{l_2}(\mathbb{R}^{N})$,
                                            in $L^{\sigma_2\widetilde{l_{\overline{\Gamma}}}}(\mathbb{R}^{N})$
                                        and in $L^{\sigma_2\widetilde{m_{\overline{\Gamma}}}}(\mathbb{R}^{N})$,
 $\bar{v}_n(x)\rightarrow \bar{v}(x)$ a.e. in $\mathbb{R}^{N}$.\\
 To get the contradiction,  we will first assume that both  $[\bar{u}\neq0]:=\{x \in \mathbb{R}^{N}: \bar u(x)\neq0\}$
and    $[\bar{v}\neq0]:=\{x \in \mathbb{R}^{N}: \bar v(x)\neq 0\}$ have zero Lebesgue measure, that is, $\bar{u}=0$ a.e.  in $\mathbb{R}^N$
 and $\bar{v}=0$ a.e.  in $\mathbb{R}^N$.
  By Lemma 2.4 and the inequality (66) in \cite{Xie2018}, we have
 \begin{align}\label{3.1.11}
 &          \frac{\min\{C_{1,3},1\}}{2^{m_{1}-1}}
            \min\left\{\|u_n\|_{1}^{l_1},\|u_n\|_{1}^{m_1}\right\}
         +  \frac{\min\{C_{2,3},1\}}{2^{m_{2}-1}}
            \min\left\{\|v_n\|_{2}^{l_2},\|v_n\|_{2}^{m_2}\right\}
            \nonumber\\
&~~~~
         -\min\{C_{1,3},1\}
         -\min\{C_{2,3},1\}
             \nonumber\\
&\leq      C_{1,3} \min\left\{ \|\nabla u_n\|_{\Phi_1}^{l_1},\|\nabla u_n\|_{\Phi_1}^{m_{1}}\right\}
        +   C_{2,3} \min\left\{ \|\nabla v_n\|_{\Phi_2}^{l_2},\|\nabla v_n\|_{\Phi_2}^{m_{2}}\right\}
            \nonumber\\
&~~~~
        +  \min\left\{ \|u_n\|_{\Phi_1,V_{1}}^{l_1},\|u_n\|_{\Phi_1,V_{1}}^{m_{1}}\right\}
        +  \min\left\{ \|v_n\|_{\Phi_2,V_{2}}^{l_2},\|v_n\|_{\Phi_2,V_{2}}^{m_{2}}\right\}
         \nonumber\\
 &\leq     C_{1,3}\int_{\mathbb{R}^{N}}\Phi_1(|\nabla u_n|)dx
         + C_{2,3}\int_{\mathbb{R}^{N}}\Phi_2(|\nabla v_n|)dx
          + \int_{\mathbb{R}^{N}}V_{1}(x)\Phi_1(|u_n|)dx
          + \int_{\mathbb{R}^{N}}V_{2}(x)\Phi_2(|v_n|)dx\nonumber\\
&\leq     \widehat{M}_{1}\left(\int_{\mathbb{R}^{N}}\Phi_1(|\nabla u_n|)dx\right)
         + \widehat{M}_{2}\left(\int_{\mathbb{R}^{N}}\Phi_2(|\nabla v_n|)dx\right)
          + \int_{\mathbb{R}^{N}}V_{1}(x)\Phi_1(|u_n|)dx
          + \int_{\mathbb{R}^{N}}V_{2}(x)\Phi_2(|v_n|)dx\nonumber\\
 &=    I(u_n,v_n)+ \int_{\mathbb{R}^{N}}F(x,u_n,v_n)dx.
 \end{align}
 When $n$ large enough, we have
 \begin{align}\label{3.1.12}
 \|u_n\|_{1}^{l_1} + \|v_n\|_{2}^{l_2}
\leq     D_{1}I(u_n,v_n)
        + D_{1}\int_{\mathbb{R}^{N}}F(x,u_n,v_n)dx
        + D_{2},
 \end{align}
 where $D_{1}=\frac{1}{ \min\left\{\frac{\min\{C_{1,3},1\}}{2^{m_{1}-1}}, \frac{\min\{C_{2,3},1\}}{2^{m_{2}-1}}\right\}}$ and
 $
 D_{2}=\frac{\min\{C_{1,3},1\} +  \min\{C_{2,3},1\}
             }
             { \min\left\{\frac{\min\{C_{1,3},1\}}{2^{m_{1}-1}}, \frac{\min\{C_{2,3},1\}}{2^{m_{2}-1}}\right\}
             }
 $.
Then
 \begin{align}\label{3.1.13}
           1
 &\leq    \frac{D_{1}I(u_n,v_n)+D_{2}}
               {  \|u_n\|_{1}^{l_1}
                 +
                  \|v_n\|_{2}^{l_2}
               }
     +        \left(\int_{|(u_n,v_n)|\leq R}+\int_{|(u_n,v_n)|> R}\right)
            \frac{D_{1}F(x,u_n,v_n)}
                 { \|u_n\|_{1}^{l_1} + \|v_n\|_{2}^{l_2}
                 }dx
\nonumber\\
 &=         o_n(1)
          + \int_{|(u_n,v_n)|\leq R}
               \frac{D_{1}F(x,u_n,v_n)}
                    {\|u_n\|_{1}^{l_1} + \|v_n\|_{2}^{l_2}
                    }dx
          + \int_{|(u_n,v_n)|> R}
                    \frac{D_{1}F(x,u_n,v_n)}
                         {\|u_n\|_{1}^{l_1} + \|v_n\|_{2}^{l_2}
                         }dx,
 \end{align}
 where $R$ is a positive constant with $R>r$.
  By $(F_2)$, there exists a  constant $0<\delta<1$ such that
 \begin{equation}\label{3.1.14-}
 \frac{|F(x,u,v)|}{\Phi_{1}(|u|)+\Phi_{2}(|v|)}< C_{3}+1, \quad \forall x\in  \mathbb{R}^N, \quad 0<|(u,v)|\leq \delta.
 \end{equation}
 By $(\phi_4)$ and (\ref{3.1.14-}), we have
 \begin{equation}\label{3.1.14----}
  \frac{|F(x,u,v)|}{|u|^{l_{1}}+|v|^{l_{2}}}=\frac{|F(x,u,v)|}{\Phi_{1}(|u|)+\Phi_{2}(|v|)}\cdot \frac{\Phi_{1}(|u|)+\Phi_{2}(|v|)}{|u|^{l_{1}}+|v|^{l_{2}}}< (C_{3}+1)\max\{c_{12},c_{22}\}, \quad \forall x\in  \mathbb{R}^N, \quad 0<|(u,v)|\leq \delta.
 \end{equation}
By the fact that $F$ and $\Phi_{i}$ are continuous,  there exist two constants $\overline{C}_{10}>0$ and $R>0$ such that
 \begin{equation}\label{3.1.14--}
 \frac{|F(x,u,v)|}{|u|^{l_{1}}+|v|^{l_{2}}}\le \frac{\max_{\delta\leq|(u,v)|\leq R}|F(x,u,v)|}{\min_{\delta\leq|(u,v)|\leq R}|u|^{l_{1}}+|v|^{l_{2}}}
                                           =:\overline{C}_{10}, \quad \forall x\in  \mathbb{R}^N, \quad \delta\leq|(u,v)|\leq R,
 \end{equation}
 combining with (\ref{3.1.14----}), which implies that there exists  a positive constant $C_{10}$ such that
\begin{equation}\label{3.1.14---}
 \frac{|F(x,u,v)|}{|u|^{l_{1}}+|v|^{l_{2}}}\le C_{10}, \quad \forall x\in  \mathbb{R}^N, \quad 0<|(u,v)|\leq R.
 \end{equation}
 Set
 \begin{eqnarray}
 &  &  \bar{B}_{n,R}=\{x\in\mathbb{R}^N||(u_n(x),v_n(x))|\leq R\},\nonumber\\
 &  &   \Omega_{1n}=\left\{x\in \bar{B}_{n,R}|u_n(x)=0 \mbox{ and } v_n(x)=0\right\},\nonumber\\
 &  &  \Omega_{2n}=\left\{x\in \bar{B}_{n,R}|u_n(x)\not=0 \mbox{ and } v_n(x)=0\right\},\nonumber\\
 &  &   \Omega_{3n}=\left\{x\in \bar{B}_{n,R}|u_n(x)=0 \mbox{ and } v_n(x)\not=0\right\},\nonumber\\
 &  &  \Omega_{4n}=\left\{x\in \bar{B}_{n,R}|u_n(x)\not=0 \mbox{ and } v_n(x)\not=0\right\}.\nonumber
 \end{eqnarray}
 Then by $(F_0)$,
 \begin{align}\label{a1}
   \int_{\Omega_{1n}}\frac{F(x,u_n,v_n)}{\|u_n\|_{1}^{l_1}+\|v_n\|_{2}^{l_2}}dx
& = 0.
 \end{align}
 Note that $v_n(x)=0$ on  $\Omega_{2n}$. We have
 \begin{align}\label{a2}
   \int_{\Omega_{2n}}\frac{F(x,u_n,v_n)}{\|u_n\|_{1}^{l_1}+\|v_n\|_{2}^{l_2}}dx
& \le  \int_{\Omega_{2n}}\frac{F(x,u_n,0)}{\|u_n\|_{1}^{l_1}}dx=\int_{\Omega_{2n}}\frac{F(x,u_n,0)}{|u_n|^{l_1}}\cdot|\bar{u}_n|^{l_1}dx\nonumber\\
&  =  \int_{\Omega_{2n}}\frac{F(x,u_n,v_n(x))}{|u_n|^{l_1}+|v_n|^{l_1}}\cdot|\bar{u}_n|^{l_1}dx\nonumber\\
& \le C_{10} \int_{\mathbb{R}^N}|\bar{u}_n|^{l_1}dx\to 0, \mbox{  as } n\to \infty.
 \end{align}
 Similarly, we also have
  \begin{align}\label{a3}
   \int_{\Omega_{3n}}\frac{F(x,u_n,v_n)}{\|u_n\|_{1}^{l_1}+\|v_n\|_{2}^{l_2}}dx \to 0, \mbox{  as } n\to \infty.
 \end{align}
 Moreover,
 \begin{align}\label{a4}
&
   \int_{\Omega_{4n}}\frac{F(x,u_n,v_n)}{\|u_n\|_{1}^{l_1}+\|v_n\|_{2}^{l_2}}dx \nonumber\\
&\qquad\qquad
   = \int_{\Omega_{4n}}
             \frac{F(x,u_n,v_n)}{\frac{|u_n|^{l_1}}{|\bar{u}_n|^{l_1}}+\frac{|v_n|^{l_2}}{|\bar{v}_n|^{l_2}}}dx\nonumber\\
&\qquad\qquad
   \leq \int_{\Omega_{4n}}
             \frac{F(x,u_n,v_n)}{\left(|u_n|^{l_1}+|v_n|^{l_2}\right)\min\left\{\frac{1}{|\bar{u}_n|^{l_1}},\frac{1}{|\bar{v}_n|^{l_2}}\right\}}dx
             \nonumber\\
&\qquad\qquad
   \leq \int_{\Omega_{4n}}
             \frac{F(x,u_n,v_n)\max\left\{|\bar{u}_n|^{l_1},|\bar{v}_n|^{l_2}\right\}}{|u_n|^{l_1}+|v_n|^{l_2}}dx
             \nonumber\\
&\qquad\qquad
   \leq \int_{\Omega_{4n}}
             \frac{F(x,u_n,v_n)}{|u_n|^{l_1}+|v_n|^{l_2}}(|\bar{u}_n|^{l_1}+|\bar{v}_n|^{l_2})dx\nonumber\\
&\qquad\qquad
   \leq C_{10}\int_{\mathbb{R}^N}
            (|\bar{u}_n|^{l_1}+|\bar{v}_n|^{l_2})dx
\rightarrow 0, \mbox{  as } n\to \infty.
 \end{align}
  Then by (\ref{a1})-(\ref{a4}), we have
 \begin{align}\label{3.1.14}
   \int_{\bar{B}_{n,R}}\frac{F(x,u_n,v_n)}{\|u_n\|_{1}^{l_1}+\|v_n\|_{2}^{l_2}}dx=\sum_{i=1}^4 \int_{\Omega_{in}}\frac{F(x,u_n,v_n)}{\|u_n\|_{1}^{l_1}+\|v_n\|_{2}^{l_2}}dx=o(1).
 \end{align}
 Next, we set
 \begin{eqnarray}
 &  &  \Omega'_{2n}=\left\{x\in \mathbb{R}^N/\bar{B}_{n,R}|u_n(x)\not=0 \mbox{ and } v_n(x)=0\right\},\nonumber\\
 &  &   \Omega'_{3n}=\left\{x\in \mathbb{R}^N/\bar{B}_{n,R}|u_n(x)=0 \mbox{ and } v_n(x)\not=0\right\},\nonumber\\
 &  &  \Omega'_{4n}=\left\{x\in \mathbb{R}^N/\bar{B}_{n,R}|u_n(x)\not=0 \mbox{ and } v_n(x)\not=0\right\}.\nonumber
 \end{eqnarray}
 Then for large $n$, we have
 \begin{align}\label{3.1.15}
      &             \int_{ \Omega'_{4n}}\frac{F(x,u_n,v_n)}{\|u_n\|_{1}^{l_1}+\|v_n\|_{2}^{l_2}}dx
      \nonumber\\
 \leq   &                    \int_{ \Omega'_{4n}}
                 \frac{|F(x,u_n,v_n)|}
                      {   \|u_n\|_{1}^{\sigma_1}
                        + \|v_n\|_{2}^{\sigma_2}
                      }
                 dx
                 \nonumber\\
 \leq &                  \int_{ \Omega'_{4n}}
                \frac{|F(x,u_n,v_n)|}{|u_n|^{\sigma_1}+|v_n|^{\sigma_2}}
                (|\bar{u}_n|^{\sigma_1}+|\bar{v}_n|^{\sigma_2})
              dx.
 \end{align}
 If
 $\frac{|F(x,u_n,v_n)|}{|u_n|^{\sigma_1}+|v_n|^{\sigma_2}} \geq 1,$
 then by H\"{o}lder's inequality,  (\ref{3.1.10}) and Remark 3.2, we have
  \begin{align}\label{3.1.15+}
& \int_{ \Omega'_{4n}}
                \frac{|F(x,u_n,v_n)|}{|u_n|^{\sigma_1}+|v_n|^{\sigma_2}}
                (|\bar{u}_n|^{\sigma_1}+|\bar{v}_n|^{\sigma_2})
              dx  \nonumber\\
&~~~~\le\left(\int_{ \Omega'_{4n}}
                \left(\frac{|F(x,u_n,v_n)|}{|u_n|^{\sigma_1}+|v_n|^{\sigma_2}}\right)^{l_{\overline{\Gamma}}}dx\right)^{\frac{1}{l_{\overline{\Gamma}}}}
                \left(\int_{ \Omega'_{4n}}(|\bar{u}_n|^{\sigma_1}+|\bar{v}_n|^{\sigma_2})^{\widetilde{l}_{\overline{\Gamma}}}dx\right)^{\frac{1}{\widetilde{l}_{\overline{\Gamma}}}}\nonumber\\
&~~~~
 \leq  2^{\frac{\widetilde{l}_{\overline{\Gamma}}-1}{\widetilde{l}_{\overline{\Gamma}}}}  \left(\frac{1}{\overline{\Gamma}(1)}\int_{ \Omega'_{4n}}
                \overline{\Gamma}\left(\frac{|F(x,u_n,v_n)|}{|u_n|^{\sigma_1}+|v_n|^{\sigma_2}}\right)
              dx\right)^{\frac{1}{l_{\overline{\Gamma}}}}\left(  \|\bar{u}_n\|^{\sigma_1}_{L^{\sigma_1\widetilde{l}_{\overline{\Gamma}}}}
                   + \|\bar{v}_n\|^{\sigma_2}_{L^{\sigma_2\widetilde{l}_{\overline{\Gamma}}}}
           \right)\nonumber\\
&~~~~
 \leq  2^{\frac{\widetilde{l}_{\overline{\Gamma}}-1}{\widetilde{l}_{\overline{\Gamma}}}}
              \left(\frac{C_4}{\overline{\Gamma}(1)}\int_{ \mathbb{R}^N}
                \overline{F}(x,u_n,v_n)
              dx\right)^{\frac{1}{l_{\overline{\Gamma}}}}
             \left(  \|\bar{u}_n\|^{\sigma_1}_{L^{\sigma_1\widetilde{l}_{\overline{\Gamma}}}}
                   + \|\bar{v}_n\|^{\sigma_2}_{L^{\sigma_2\widetilde{l}_{\overline{\Gamma}}}}
           \right)\nonumber\\
 &~~~~
 \leq  2^{\frac{\widetilde{l}_{\overline{\Gamma}}-1}{\widetilde{l}_{\overline{\Gamma}}}}
              \left(\frac{C_4(c+1)}{\overline{\Gamma}(1)}\right)^{\frac{1}{l_{\overline{\Gamma}}}}
             \left(  \|\bar{u}_n\|^{\sigma_1}_{L^{\sigma_1\widetilde{l}_{\overline{\Gamma}}}}
                   + \|\bar{v}_n\|^{\sigma_2}_{L^{\sigma_2\widetilde{l}_{\overline{\Gamma}}}}
           \right)
\rightarrow 0, \mbox{ as } n\to \infty.
 \end{align}
 Similarly, if
 $\frac{|F(x,u_n,v_n)|}{|u_n|^{\sigma_1}+|v_n|^{\sigma_2}} < 1,$
 then, by H\"{o}lder's inequality,  (\ref{3.1.10}) and Remark 3.2, we get that
 \begin{align}\label{3.1.15-}
& \int_{ \Omega'_{4n}}
                \frac{|F(x,u_n,v_n)|}{|u_n|^{\sigma_1}+|v_n|^{\sigma_2}}
                (|\bar{u}_n|^{\sigma_1}+|\bar{v}_n|^{\sigma_2})
              dx \nonumber\\
&~~~~
 \leq              2^{\frac{\widetilde{m}_{\overline{\Gamma}}-1}{\widetilde{m}_{\overline{\Gamma}}}}
              \left(\frac{C_4}{\overline{\Gamma}(1)}\int_{ \mathbb{R}^N}
                \overline{F}(x,u_n,v_n)
              dx\right)^{\frac{1}{m_{\overline{\Gamma}}}}
              \left(  \|\bar{u}_n\|^{\sigma_1}_{L^{\sigma_1\widetilde{m}_{\overline{\Gamma}}}}
                   + \|\bar{v}_n\|^{\sigma_2}_{L^{\sigma_2\widetilde{m}_{\overline{\Gamma}}}}
           \right)\nonumber\\
 &~~~~
 \leq              2^{\frac{\widetilde{m}_{\overline{\Gamma}}-1}{\widetilde{m}_{\overline{\Gamma}}}}
              \left(\frac{C_4(c+1)}{\overline{\Gamma}(1)}\right)^{\frac{1}{m_{\overline{\Gamma}}}}
              \left(  \|\bar{u}_n\|^{\sigma_1}_{L^{\sigma_1\widetilde{m}_{\overline{\Gamma}}}}
                   + \|\bar{v}_n\|^{\sigma_2}_{L^{\sigma_2\widetilde{m}_{\overline{\Gamma}}}}
           \right)
\rightarrow 0, \mbox{ as } n\to \infty.
\end{align}
Hence, (\ref{3.1.15})-(\ref{3.1.15-})  imply that
 \begin{align}\label{3.1.15x}
 \int_{\Omega'_{4n}}\frac{F(x,u_n,v_n)}{\|u_n\|_{1}^{l_1}+\|v_n\|_{2}^{l_2}}dx
 \to 0, \mbox{ as } n\to \infty.
 \end{align}
Note that $v_n(x)=0$ on $ \Omega'_{2n}$. Then we  have
 \begin{align}\label{b1}
      &             \int_{ \Omega'_{2n}}\frac{F(x,u_n,v_n)}{\|u_n\|_{1}^{l_1}+\|v_n\|_{2}^{l_2}}dx
      \nonumber\\
 \leq   &                    \int_{ \Omega'_{2n}}
                 \frac{|F(x,u_n,0)|}
                      {   \|u_n\|_{1}^{\sigma_1}
                        + \|v_n\|_{2}^{\sigma_2}
                      }
                 dx
                 \nonumber\\
 \leq   &                    \int_{ \Omega'_{2n}}
                 \frac{|F(x,u_n,0)|}
                      {   \|u_n\|_{1}^{\sigma_1}
                                            }
                 dx
                 \nonumber\\
= &                    \int_{ \Omega'_{2n}}
                \frac{|F(x,u_n,v_n)|}{|u_n|^{\sigma_1}+|v_n|^{\sigma_2}}\cdot
                |\bar{u}_n|^{\sigma_1}
              dx.
 \end{align}
 Similar to the arguments of (\ref{3.1.15+}) and (\ref{3.1.15-}), we can also get
 \begin{align}
 \int_{ \Omega'_{2n}}
                \frac{|F(x,u_n,v_n)|}{|u_n|^{\sigma_1}+|v_n|^{\sigma_2}}\cdot
                |\bar{u}_n|^{\sigma_1}
              dx\to 0, \mbox{ as } n\to \infty,\nonumber
               \end{align}
which shows that
 \begin{align}\label{b2}
 \int_{ \Omega'_{2n}}\frac{F(x,u_n,v_n)}{\|u_n\|_{1}^{l_1}+\|v_n\|_{2}^{l_2}}dx\to 0, \mbox{ as } n\to \infty.
  \end{align}
  Similarly, we also get
  \begin{align}\label{b3}
 \int_{ \Omega'_{3n}}\frac{F(x,u_n,v_n)}{\|u_n\|_{1}^{l_1}+\|v_n\|_{2}^{l_2}}dx\to 0, \mbox{ as } n\to \infty.
  \end{align}
Hence, (\ref{3.1.15x}), (\ref{b2}) and (\ref{b3})  imply that
 \begin{align}\label{b4}
 \int_{\mathbb{R}^N/\bar{B}_{n,R}}\frac{F(x,u_n,v_n)}{\|u_n\|_{1}^{l_1}+\|v_n\|_{2}^{l_2}}dx
 =o(1).
 \end{align}
By combining \eqref{3.1.14}, \eqref{b4} with \eqref{3.1.13}, we get a contradiction.
\par
Next, we  assume that
      $[\bar{u}\neq0]$
or    $[\bar{v}\neq0]$ has nonzero Lebesgue measure.
It is clear that
 $$
 |u_n|=|\bar{u}_n|\|u_n\|_{1}\rightarrow +\infty \quad \mbox{ for\;all }[\bar{u}\neq 0]
 $$
 or
 $$
 |v_n|=|\bar{v}_n|\|v_n\|_{2}\rightarrow +\infty \quad \mbox{  for\;all }[\bar{v}\neq 0].
 $$
Moreover,  by (2) in Lemma 2.3, $F(x,u,v)\geq 0$ and the assumption $(F_3)$  show that
 \begin{align}\label{3.1.15---}
& \lim_{|(u,v)|\rightarrow+\infty}\frac{F(x,u,v)}{|u|^{\sigma_1}+|v|^{\sigma_2}} \nonumber\\
&~~~~
 = \begin{cases}
              \lim_{|(u,v)|\rightarrow+\infty}\frac{F(x,u,v)}{|u|^{\sigma_1}+|v|^{\sigma_2}},
                   \;\;\; \;\;\; \;\; \;             \text { if }\;\;|u|\geq 1, |v|\geq 1, |(u,v)|\rightarrow+\infty, \\
              \lim_{|(u,v)|\rightarrow+\infty}\frac{F(x,u,v)}{|u|^{\sigma_1}+|v|^{\sigma_2}}
         \geq \lim_{|(u,v)|\rightarrow+\infty}\frac{F(x,u,v)}{1+|v|^{\sigma_2}},
                   \;\;\; \;\;\; \;\; \;             \text { if }\;\;0\leq |u|< 1, |v|\rightarrow+\infty,\\
              \lim_{|(u,v)|\rightarrow+\infty}\frac{F(x,u,v)}{|u|^{\sigma_1}+|v|^{\sigma_2}}
         \geq  \lim_{|(u,v)|\rightarrow+\infty}\frac{F(x,u,v)}{1+|u|^{\sigma_1}},
                   \;\;\; \;\;\; \;\; \;             \text { if }\;\;0\leq |v|< 1, |u|\rightarrow+\infty
      \end{cases} \nonumber\\
&~~~~
 \geq \begin{cases}
              \lim_{|(u,v)|\rightarrow+\infty}\frac{F(x,u,v)}{|u|^{m_1}+|v|^{m_2}},
                   \;\;\; \;\;\; \;\; \;             \text { if }\;\;|u|\geq 1, |v|\geq 1, \\
                \lim_{|(u,v)|\rightarrow+\infty}\frac{F(x,u,v)}{|v|^{m_2}+|v|^{m_2}}
        \geq  \frac{1}{2}\lim_{|(u,v)|\rightarrow+\infty}\frac{F(x,u,v)}{|u|^{m_1}+|v|^{m_2}},
                   \;\;\; \;\;\; \;\; \;             \text { if }\;\;0\leq |u|< 1, |v|>1,\\
             \lim_{|(u,v)|\rightarrow+\infty}\frac{F(x,u,v)}{|u|^{m_1}+|u|^{m_1}}
        \geq  \frac{1}{2}\lim_{|(u,v)|\rightarrow+\infty}\frac{F(x,u,v)}{|u|^{m_1}+|v|^{m_2}},
                   \;\;\; \;\;\; \;\; \;             \text { if }\;\;0\leq |v|< 1, |u|>1
      \end{cases}  \nonumber\\
&~~~~
= +\infty \quad \mbox{  for \; a.e. } x\in G.
\end{align}
Hence, combining with $(F_4)$, it is easy to see that
\begin{align}\label{3.1.15-----}
\lim_{|(u,v)|\rightarrow +\infty}\overline{F}(x,u,v)
&~~~~
 \geq \frac{1}{C_{4}} \lim_{|(u,v)|\rightarrow +\infty}\overline{\Gamma}\left(\frac{F(x,u,v)}{|u|^{\sigma_1}+|v|^{\sigma_2}}\right)
 \nonumber\\
&~~~~
 \geq \frac{1}{C_{4}} \lim_{|(u,v)|\rightarrow +\infty}\overline{\Gamma}(1)\max\left\{\left(\frac{F(x,u,v)}{|u|^{\sigma_1}+|v|^{\sigma_2}}\right)^{l_{\overline{\Gamma}}},
                                        \left(\frac{F(x,u,v)}{|u|^{\sigma_1}+|v|^{\sigma_2}}\right)^{m_{\overline{\Gamma}}}
                                 \right\}
\nonumber\\
&~~~~
 \geq \frac{1}{C_{4}} \overline{\Gamma}(1)
     \max\left\{\left(\lim_{|(u,v)|\rightarrow+\infty}\frac{F(x,u,v)}{|u|^{\sigma_1}+|v|^{\sigma_2}}\right)^{l_{\overline{\Gamma}}},
                \left(\lim_{|(u,v)|\rightarrow +\infty}\frac{F(x,u,v)}{|u|^{\sigma_1}+|v|^{\sigma_2}}\right)^{m_{\overline{\Gamma}}}
                                 \right\}\nonumber\\
&~~~~
 =+\infty \quad \mbox{  for \; a.e. } x\in G.
 \end{align}
 Then, by \eqref{3.1.10},  Fatou Lemma and the above formula, we have
 $$
      c+1
 \geq  \lim_{n\rightarrow +\infty}\int_{\mathbb{R}^N} \overline{F}(x,u_n,v_n)dx
 \geq  \lim_{n\rightarrow +\infty}\int_{G} \overline{F}(x,u_n,v_n)dx
 \geq    \int_{G} \lim_{n\rightarrow +\infty}\overline{F}(x,u_n,v_n)dx
 =     +\infty,$$
 which is a contradiction.
 So, both $\|u_n\|_{1}\rightarrow +\infty$ and
 $\|v_n\|_{2}\rightarrow +\infty$ do not hold.

\vskip2mm
\noindent {\bf Case 2}.
Suppose that $\|u_n\|_{1}\leq D_{3}$
or
             $\|v_n\|_{2}\leq D_{3}$
for some $D_{3}>0$ and all $n\in\mathbb{N}$.
Without loss of generality, we assume that
    $\|u_n\|_{1}\rightarrow+\infty$
and
    $\|v_n\|_{2}\leq D_{3}$
for some $D_{3}>0$ and all $n\in\mathbb{N}$.
Let $\bar{u}_n=\frac{u_n}{\|u_n\|_{1}}$ and $\bar{v}_n=\frac{v_n}{\|u_n\|_{1}}$.
Then $\|{\bar{u}_n}\|_{1}=1$ and $\|{\bar{v}_n}\|_{2}\rightarrow 0$.
Passing to a subsequences $\{(\bar{u}_n,\bar{v}_n)\}$, by Remark 3.2,
there exist $\bar{u}\in W_{1}$ and $v\in W_{2}$ such that\\
$\star$
\quad  $\bar{u}_n\rightharpoonup \bar{u}$    in $W_{1}$,
\quad      $\bar{u}_n\rightarrow \bar{u}$    in $L^{l_1}(\mathbb{R}^{N})$,
                                             in $L^{\sigma_1\widetilde{l_{\overline{\Gamma}}}}(\mathbb{R}^{N})$
                                      and in $L^{\sigma_1\widetilde{m_{\overline{\Gamma}}}}(\mathbb{R}^{N})$,
\quad      $\bar{u}_n(x)\rightarrow \bar{u}(x)$  a.e. in $\mathbb{R}^{N}$;\\
$\star$
\quad  $\bar{v}_n\rightarrow0$       in $W_{2}$,
\quad  $\bar{v}_n\rightarrow 0$      in $L^{l_2}(\mathbb{R}^{N})$,
                                      in $L^{\sigma_2\widetilde{l_{\overline{\Gamma}}}}(\mathbb{R}^{N})$
                                        and in $L^{\sigma_2\widetilde{m_{\overline{\Gamma}}}}(\mathbb{R}^{N})$,
\quad    $\bar{v}_n(x)\rightarrow 0$   a.e. in $\mathbb{R}^{N}$;\\
$\star$
\quad  $v_n\rightharpoonup v$     in $W_{2}$,
\quad  $v_n\rightarrow v$      in $L^{l_2}(\mathbb{R}^{N})$,
                               in $L^{\sigma_2\widetilde{l_{\overline{\Gamma}}}}(\mathbb{R}^{N})$
                                        and in $L^{\sigma_2\widetilde{m_{\overline{\Gamma}}}}(\mathbb{R}^{N})$,
\quad  $v_n(x)\rightarrow v(x)$   a.e. in $\mathbb{R}^{N}$.\\
 Firstly, we assume that $[\bar{u}\neq0]$ has nonzero Lebesgue measure. We can see that
 $$|u_n|=|\bar{u}_n|\|u_n\|_{1}\rightarrow +\infty \quad \mbox{ for\;all }[\bar{u}\neq 0].$$
 Then, being analogue to Case 1, we get a contradiction by
 $$c+1\geq\int_{\mathbb{R}^N} \overline{F}(x,u_n,v_n)dx\rightarrow +\infty.$$
 Next, we suppose that $[\bar{u}\neq0]$ has zero Lebesgue measure, that is, $\bar{u}=0$ a.e. in $\mathbb{R}^N$.
 By \eqref{3.1.10}, we can see that
  \begin{align}\label{3.1.10+}
\int_{\mathbb{R}^N}\overline{F}(x,u_n,v_n)dx \leq  c+1.
\end{align}
 Then when $n$ large enough, we can choose two positive constants
$D_4,D_5$
 such that \eqref{3.1.11} is changed into
 \begin{align}\label{d5}
 \|u_n\|_{1}^{l_1} \leq D_4I(u_n,v_n)+D_4\int_{\mathbb{R}^N}F(x,u_n,v_n)dx+D_{5}.
 \end{align}
 Similar to the arguments of (\ref{a1})-(\ref{a3}),  for any given
  $$
 D_6>\max\left\{5C_{10}C_{2,5}^{l_2}D_3^{l_2},
 5\left(\frac{C_4(c+1)}{\overline{\Gamma}(1)}\right)^{\frac{1}{l_{\overline{\Gamma}}}}C_{2,5}^{\sigma_2}D_3^{\sigma_2},
 5\left(\frac{C_4(c+1)}{\overline{\Gamma}(1)}\right)^{\frac{1}{m_{\overline{\Gamma}}}}C_{2,5}^{\sigma_2}D_3^{\sigma_2}\right\},
 $$
  we obtain that
  \begin{align}
  \label{b5} & \int_{\Omega_{1n}}\frac{F(x,u_n,v_n)}{\|u_n\|_{1}^{l_1}+D_6}dx\le \int_{\Omega_{1n}}\frac{F(x,u_n,v_n)}{\|u_n\|_{1}^{l_1}}dx =0,\\
   \label{b6}& \int_{\Omega_{2n}}\frac{F(x,u_n,v_n)}{\|u_n\|_{1}^{l_1}+D_6}dx\le\int_{\Omega_{2n}}\frac{F(x,u_n,v_n)}{\|u_n\|_{1}^{l_1}}dx \to 0, \mbox{  as } n\to \infty
   \end{align}
 and
  \begin{align}
  \label{b7} \int_{\Omega_{3n}}\frac{F(x,u_n,v_n)}{\|u_n\|_{1}^{l_1}+D_6}dx
  &\le\int_{\Omega_{3n}}\frac{F(x,u_n,v_n)}{|{v}_n|^{l_2}}\cdot\frac{|{v}_n|^{l_2}}{D_6}dx
  =\int_{\Omega_{3n}}\frac{F(x,u_n,v_n)}{|{u}_n|^{l_1}+|{v}_n|^{l_2}}\cdot\frac{|{v}_n|^{l_2}}{D_6}dx \nonumber\\
  &\leq\frac{C_{10}}{D_6}\int_{\mathbb{R}^N}|{v}_n|^{l_2}dx
  \leq\frac{C_{10}C_{2,5}^{l_2}D_3^{l_2}}{D_6}
  <\frac{1}{5}, \mbox{  as } n\to \infty.
 \end{align}
 Moreover, we  have
 \begin{align}\label{b8}
&
   \int_{\Omega_{4n}}\frac{F(x,u_n,v_n)}{\|u_n\|_{1}^{l_1}+D_6}dx\nonumber\\
&\qquad\qquad
   \leq \int_{\Omega_{4n}}
             \frac{F(x,u_n,v_n)}{\frac{|u_n|^{l_1}}{|\bar{u}_n|^{l_1}}+D_6\frac{|v_n|^{l_2}}{|{v}_n|^{l_2}}}dx\nonumber\\
&\qquad\qquad
   \leq \int_{\Omega_{4n}}
             \frac{F(x,u_n,v_n)}{\left(|u_n|^{l_1}+|v_n|^{l_2}\right)\min\left\{\frac{1}{|\bar{u}_n|^{l_1}},\frac{D_6}{|{v}_n|^{l_2}}\right\}}dx
             \nonumber\\
&\qquad\qquad
   \leq \int_{\Omega_{4n}}
             \frac{F(x,u_n,v_n)\max\left\{|\bar{u}_n|^{l_1},\frac{|{v}_n|^{l_2}}{D_6}\right\}}{|u_n|^{l_1}+|v_n|^{l_2}}dx
             \nonumber\\
&\qquad\qquad
   \leq \int_{\Omega_{4n}}
             \frac{F(x,u_n,v_n)}{|u_n|^{l_1}+|v_n|^{l_2}}(|\bar{u}_n|^{l_1}+\frac{|{v}_n|^{l_2}}{D_6})dx\nonumber\\
&\qquad\qquad
   \leq C_{10}\int_{\mathbb{R}^N}
            |\bar{u}_n|^{l_1}dx+\frac{C_{10}}{D_6}\int_{\mathbb{R}^N}|{v}_n|^{l_2}dx\nonumber\\
&\qquad\qquad
   \leq C_{10}\int_{\mathbb{R}^N}
            |\bar{u}_n|^{l_1}dx+\frac{C_{10}C_{2,5}^{l_2}D_3^{l_2}}{D_6}\nonumber\\
&\qquad\qquad
  \leq o(1)+\frac{1}{5}.
 \end{align}
 Hence, (\ref{b5})-(\ref{b8}) imply that
 \begin{align}\label{c1}
   \int_{\bar{B}_{n,R}}\frac{F(x,u_n,v_n)}{\|u_n\|_{1}^{l_1}+D_6}dx=\sum_{i=1}^4 \int_{\Omega_{in}}\frac{F(x,u_n,v_n)}{\|u_n\|_{1}^{l_1}+D_6}dx=o(1)+\frac{2}{5}.
 \end{align}
 Moreover,
 \begin{align}\label{d4}
  &           \int_{\Omega'_{4n}}\frac{F(x,u_n,v_n)}{\|u_n\|_{1}^{l_1}+D_6}dx\nonumber\\
  &\leq     \int_{\Omega'_{4n}}\frac{F(x,u_n,v_n)}{\|u_n\|_{1}^{\sigma_1}+D_6}dx\nonumber\\
  &=       \int_{\Omega'_{4n}}
         \frac{F(x,u_n,v_n)}{\frac{|u_n|^{\sigma_1}}{|\bar{u}_n|^{\sigma_1}}+\frac{D_6|v_n|^{\sigma_2}}{|v_n|^{\sigma_2}}}dx\nonumber\\
  &\leq      \int_{\Omega'_{4n}}
         \frac{F(x,u_n,v_n)}{(|u_n|^{\sigma_1}+|v_n|^{\sigma_2})\min\left\{\frac{1}{|\bar{u}_n|^{\sigma_1}},\frac{D_6}{|v_n|^{\sigma_2}}\right\}}dx\nonumber\\
  &=     \int_{\Omega'_{4n}}
         \frac{F(x,u_n,v_n)}{(|u_n|^{\sigma_1}+|v_n|^{\sigma_2})}\max\left\{|\bar{u}_n|^{\sigma_1},\frac{|v_n|^{\sigma_2}}{D_6}\right\}dx\nonumber\\
  &\leq
            \int_{\Omega'_{4n}}
             \frac{F(x,u_n,v_n)}{|u_n|^{\sigma_1}+|v_n|^{\sigma_2}}
             \left(|\bar{u}_n|^{\sigma_1}+\frac{|v_n|^{\sigma_2}}{D_6}\right)
             dx\nonumber\\
  &=
            \int_{\Omega'_{4n}}
             \frac{F(x,u_n,v_n)}{|u_n|^{\sigma_1}+|v_n|^{\sigma_2}}
             |\bar{u}_n|^{\sigma_1}
             dx
         +  \frac{1}{D_6}
            \int_{\Omega'_{4n}}
             \frac{F(x,u_n,v_n)}{|u_n|^{\sigma_1}+|v_n|^{\sigma_2}}
             |v_n|^{\sigma_2}
             dx\nonumber\\
  &\le
              \max\left\{\left(\frac{C_4(c+1)}{\overline{\Gamma}(1)}\right)^{\frac{1}{l_{\overline{\Gamma}}}}
              \|\bar{u}_n\|^{\sigma_1}_{L^{\sigma_1\widetilde{l}_{\overline{\Gamma}}}},\left(\frac{C_4(c+1)}{\overline{\Gamma}(1)}\right)^{\frac{1}{m_{\overline{\Gamma}}}}
              \|\bar{u}_n\|^{\sigma_1}_{L^{\sigma_1\widetilde{m}_{\overline{\Gamma}}}}\right\}\nonumber\\
  &              ~~~   + \frac{1}{D_6} \max\left\{\left(\frac{C_4(c+1)}{\overline{\Gamma}(1)}\right)^{\frac{1}{l_{\overline{\Gamma}}}}
              \|{v}_n\|^{\sigma_2}_{L^{\sigma_2\widetilde{l}_{\overline{\Gamma}}}},\left(\frac{C_4(c+1)}{\overline{\Gamma}(1)}\right)^{\frac{1}{m_{\overline{\Gamma}}}}
              \|{v}_n\|^{\sigma_2}_{L^{\sigma_2\widetilde{m}_{\overline{\Gamma}}}}\right\}
            \nonumber\\
 &\le    o(1)+\frac{1}{5}.
 \end{align}
 Note that $v_n(x)=0$ on $ \Omega'_{2n}$.  Similar to the arguments of (\ref{3.1.15+}) and (\ref{3.1.15-}), we can also get
 \begin{align}\label{d1}
      &             \int_{ \Omega'_{2n}}\frac{F(x,u_n,v_n)}{\|u_n\|_{1}^{l_1}+D_6}dx
      \nonumber\\
 \leq   &                    \int_{ \Omega'_{2n}}
                 \frac{|F(x,u_n,0)|}
                      {   \|u_n\|_{1}^{\sigma_1}
                      }
                 dx
                 \nonumber\\
 \leq &                    \int_{ \Omega'_{2n}}
                \frac{|F(x,u_n,v_n)|}{|u_n|^{\sigma_1}+|v_n|^{\sigma_2}}\cdot
                |\bar{u}_n|^{\sigma_1}
              dx\to 0, \mbox{ as } n\to \infty.
               \end{align}
  Note that $\bar{v}_n(x)=\frac{v_n(x)}{\|u_n\|_{1}}$ and ${u}_n(x)=0$ on $\Omega'_{3n}$.  Similar to the argument of (\ref{d4}), we have
  \begin{align}\label{d3}
                   \int_{ \Omega'_{3n}}\frac{F(x,u_n,v_n)}{\|u_n\|_{1}^{l_1}+D_6}dx
 \leq                       \int_{ \Omega'_{3n}}
                 \frac{|F(x,0,v_n)|}
                      {  D_6\frac{|v_n(x)|^{\sigma_2}}{|v_n(x)|^{\sigma_2}}
                      }
                 dx
 =                    \frac{1}{D_6}\int_{ \Omega'_{3n}}
                \frac{|F(x,u_n,v_n)|}{|u_n|^{\sigma_1}+|v_n|^{\sigma_2}}\cdot
                |{v}_n|^{\sigma_2}
              dx\le\frac{1}{5}.
 \end{align}
 Hence, (\ref{d4}), (\ref{d1})  and (\ref{d3})  imply that
 \begin{align}\label{d8}
 \int_{\mathbb{R}^N/\bar{B}_{n,R}}\frac{F(x,u_n,v_n)}{\|u_n\|_{1}^{l_1}+D_6}dx
 \le o(1)+\frac{2}{5}.
 \end{align}
 Then, by \eqref{d5}, \eqref{c1} and \eqref{d8}, we have
 \begin{align*}
             1
 &\leq      \frac{D_4I(u_n,v_n)+D_5+D_6}{\|u_n\|_{1}^{l_1}+D_6}
         +  D_4\int_{\mathbb{R}^N}\frac{F(x,u_n,v_n)}{\|u_n\|_{1}^{l_1}+D_6}dx\nonumber\\
  &  =      o(1)
         +  \int_{\bar{B}_{n,R}}\frac{F(x,u_n,v_n)}{\|u_n\|_{1}^{l_1}+D_{6}}dx
         +  \int_{\mathbb{R}^N/\bar{B}_{n,R}}\frac{F(x,u_n,v_n)}{\|u_n\|_{1}^{l_1}+D_{6}}dx\nonumber\\
 &\le    o(1)+\frac{4}{5},
 \end{align*}
 which is a contradiction. Based on these advantages,  we could get the conclusion of boundedness  for sequence $\{(u_n,v_n)\}$.\qed
 \vskip2mm
 \noindent
{\bf Lemma 3.8.}
 Suppose that $(\phi_1)$--$(\phi_4)$, $(V_0)$, $(V_1)$, $(M_1)$,  $(F_1)$, $(F_3)$ and $(F_4)$ hold. Then $I$ satisfies the (C)-condition.
 \vskip2mm
\noindent
 {\bf Proof.} Let $\{(u_n,v_n)\}$ be any (C)-sequence of $I$ in $W$. Lemma 3.7 shows that $\{(u_n,v_n)\}$ is bounded. Passing to a subsequence $\{(u_n,v_n)\}$, by Remark 3.2, there exists a point $(u,v)\in W$ such that\\
 $\star$\quad $u_n\rightharpoonup u$ in $W_{1}$,
 \quad $u_n\rightarrow u$ in $L^{\Psi_1}(\mathbb{R}^{N})$,
                          in $L^{\Phi_1}(\mathbb{R}^{N})$,
                          in $L^{m_1}(\mathbb{R}^{N})$,
 \quad $u_n(x)\rightarrow u(x)$ a.e. in $\mathbb{R}^{N}$;\\
 $\star$\quad $v_n\rightharpoonup v$ in $W_{2}$,
 \quad $v_n\rightarrow v$ in $L^{\Psi_2}(\mathbb{R}^{N})$,
                          in $L^{\Phi_2}(\mathbb{R}^{N})$,
                          in $L^{m_2}(\mathbb{R}^{N})$,
 \quad $v_n(x)\rightarrow v(x)$ a.e. in $\mathbb{R}^{N}$.
 \par
 Now, we define the operators
 $\mathcal{F}: W_{1}\rightarrow(W_{1})^*$ by
 $$
 \langle\mathcal{F}(u),\tilde{u}\rangle
 :=M_{1}\left(\int_{\mathbb{R}^N}\Phi_1(|\nabla u|)dx\right)
        \int_{\mathbb{R}^N}\phi_1(|\nabla u|)\nabla u\nabla \tilde{u}dx
  + \int_{\mathbb{R}^N}V_{1}(x)\phi_1(| u|) u \tilde{u}dx, \;\;
 u, \tilde{u} \in W_{1}
 $$
 and $\mathcal{G}: W_{2}\rightarrow(W_{2})^*$
 by
 $$
     \langle\mathcal{G}(v),\tilde{v}\rangle
 := M_{2}\left(\int_{\mathbb{R}^N}\Phi_2(|\nabla v|)dx\right)
            \int_{\mathbb{R}^N}\phi_2(|\nabla v|)\nabla v\nabla \tilde{v}dx
    + \int_{\mathbb{R}^N}V_{2}(x)\phi_2(| v|) v \tilde{v}dx,
            \;\;v, \tilde{v} \in W_{2}.
 $$
 Then, we have
 \begin{align}\label{3.1.19}
    \langle\mathcal{F}(u_n),u_n-u\rangle
 & =   M_{1}\left(\int_{\mathbb{R}^N}\Phi_1(|\nabla u_n|)dx\right)
        \int_{\mathbb{R}^N}\phi_1(|\nabla u_n|)\nabla u_n\nabla (u_n-u)dx
        \nonumber \\
 &~~~~~~
  + \int_{\mathbb{R}^N}V_{1}(x)\phi_1(|u_n|) u_n (u_n-u)dx\\
 &=   \langle I'(u_n,v_n),(u_n-u,0)\rangle
   +  \int_{\mathbb{R}^N}F_{u}(x,u_n,v_n)(u_n-u)dx.\nonumber
 \end{align}
 \eqref{2.2.1} and the boundness of $\{u_n\}$ show that
 \begin{equation}\label{3.1.20}
 |\langle I'(u_n,v_n),(u_n-u,0)\rangle|\leq\|I'(u_n,v_n)\|_{W^{*}}\|u_n-u\|_{1}\rightarrow 0.
 \end{equation}
  By $(F_1)$ and H\"{o}lder's inequality, we get
 \begin{align}\label{3.1.21}
 &\left|\int_{\mathbb{R}^N}F_{u}(x,u_n,v_n)(u_n-u)dx\right|\nonumber\\
 &\;\;\leq  C_2\int_{\mathbb{R}^N}(|u_n|^{l_{1}-1}+\psi_1(|u_n|)+\widetilde{\Psi}_1^{-1}(\Psi_2(|v_n|)))|u_n-u|dx
        \nonumber\\
 &\;\; =    C_2\int_{\mathbb{R}^N}|u_n|^{l_{1}-1}|u_n-u|dx
 +
        C_2\int_{\mathbb{R}^N}(\psi_1(|u_n|)+\widetilde{\Psi}_1^{-1}(\Psi_2(|v_n|)))|u_n-u|dx
        \nonumber\\
 &\;\;\leq    2C_2\||u_n|^{l_{1}-1}\|_{\widetilde{\Phi}_1}\|u_n-u\|_{\Phi_1}
        + 2C_2\|\psi_1(|u_n|)+\widetilde{\Psi}_1^{-1}(\Psi_2(|v_n|))\|_{\widetilde{\Psi}_1}\|u_n-u\|_{\Psi_1}.
 \end{align}
 Condition $(F_1)$ shows that functions $\Psi_1$ and $\widetilde{\Psi}_1$ are $N$-functions satisfying $\Delta_2$-condition globally, which together with the convexity of $N$-function, $(\phi_4)$, Lemma 2.4, Remark 2.8, Remark 3.2, inequality (A.9) in \cite{Fukagai2006}, inequality \eqref{2.1.2} and the boundedness of $\{(u_n,v_n)\}$, imply that  the  boundedness of the following integrals
  \begin{align}\label{3.1.22}
         \int_{\mathbb{R}^{N}}\widetilde{\Phi}_1(|u_n|^{l_{1}-1})dx
&  =      \int_{|u_{n}|= 0}\widetilde{\Phi}_1(|u_n|^{l_{1}-1})dx
       + \int_{0<|u_{n}|< 1}\widetilde{\Phi}_1(|u_n|^{l_{1}-1})dx
       + \int_{|u_{n}|\geq 1}\widetilde{\Phi}_1(|u_n|^{l_{1}-1})dx\nonumber\\
\leq &
         \int_{\mathbb{R}^{N}}\widetilde{\Phi}_1\left(\frac{\Phi_{1}(|u_n|)}{c_{1,1}|u_n|}\right)dx
       + \widetilde{\Phi}_1(1)\int_{\mathbb{R}^{N}}|u_n|^{l_{1}}dx
\nonumber\\
\leq &
         \int_{\mathbb{R}^{N}}\left(\frac{1}{c_{1,1}}+K_{\frac{1}{c_{1,1}}}\right)\Phi_{1}(|u_n|)dx
       + \widetilde{\Phi}_1(1)\int_{\mathbb{R}^{N}}|u_n|^{m_{1}}dx
\nonumber\\
\leq &
         \left(\frac{1}{c_{1,1}}+K_{\frac{1}{c_{1,1}}}\right)\left(\|u_n\|_{\Phi_{1}}^{l_{1}}+\|u_n\|_{\Phi_{1}}^{m_{1}}\right)
       +  \widetilde{\Phi}_1(1)\|u_n\|_{L^{m_{1}}}^{m_{1}}
 \end{align}
 and
\begin{align}\label{3.1.22-}
&        \int_{\mathbb{R}^N}\widetilde{\Psi}_1(\psi_1(|u_n|)+\widetilde{\Psi}_1^{-1}(\Psi_2(|v_n|)))dx\nonumber\\
&\qquad\qquad
 \leq    \widetilde{K}_{2}\int_{\mathbb{R}^N}\widetilde{\Psi}_1\left(\frac{\psi_1(|u_n|)}{2}+\frac{\widetilde{\Psi}_1^{-1}(\Psi_2(|v_n|))}{2}\right)dx\nonumber\\
&\qquad\qquad
 \leq    \frac{\widetilde{K}_{2}}{2}\int_{\mathbb{R}^N}\left(\widetilde{\Psi}_1\left(\psi_1(|u_n|)\right)+\widetilde{\Psi}_1\left(\widetilde{\Psi}_1^{-1}(\Psi_2(|v_n|))\right)\right)dx
 \nonumber\\
&\qquad\qquad
 \leq    \frac{\widetilde{K}_{2}}{2}\int_{\mathbb{R}^N}\left(\Psi_1\left(2|u_n|\right)+\Psi_2(|v_n|)\right)dx
 \nonumber\\
&\qquad\qquad
 \leq    \frac{\widetilde{K}_{2}}{2}\int_{\mathbb{R}^N}\left(K_{2}\Psi_1\left(|u_n|\right)+\Psi_2(|v_n|)\right)dx
 \nonumber\\
&\qquad\qquad
 \leq   D_{5}
        \int_{\mathbb{R}^N}
                 \left(
                \Psi_1\left(|u_n|\right)
              + \Psi_2(|v_n|)
                  \right)
           dx
 \nonumber\\
&\qquad\qquad
 \leq   D_{5}
 \left(\|u_n\|_{\Psi_{1}}^{l_{\Psi_{1}}}+\|u_n\|_{\Psi_{1}}^{m_{\Psi_{1}}}+\|v_n\|_{\Psi_{2}}^{l_{\Psi_{2}}}+\|v_n\|_{\Psi_{2}}^{m_{\Psi_{2}}}\right),
\end{align}
where $D_{5}=\frac{\widetilde{K}_{2}}{2}\max\left\{K_{2},1\right\},\; \widetilde{K}_{2}>0$ and  $K_{2}>0$, which shows that
 \begin{equation}\label{3.1.23}
\||u_n|^{l_{1}-1}\|_{\widetilde{\Phi}_1}\leq D_{6}
 \end{equation}
and
 \begin{equation}\label{3.1.23-}
 \|\psi_1(|u_n|)+\widetilde{\Psi}_1^{-1}(\Psi_2(|v_n|))\|_{\widetilde{\Psi}_1}\leq D_{7}
 \end{equation}
 for some $D_{6}, D_{7}>0$. Moreover, $\star$ shows that
 \begin{equation}\label{3.1.24}
 \|u_n-u\|_{\Phi_1}\rightarrow 0\;\;\mbox{and}\;\;\|u_n-u\|_{\Psi_1}\rightarrow 0.
 \end{equation}
 Then, combining \eqref{3.1.20}, \eqref{3.1.21}, \eqref{3.1.23}-\eqref{3.1.24} with \eqref{3.1.19}, we obtain
 $$\langle\mathcal{F}(u_n),u_n-u\rangle\rightarrow 0 \quad \mbox{ as } n\rightarrow \infty.$$
 By \cite[Proposition A.3]{Carvalho2015}, $\mathcal{F}$ is of the class $(S_+)$, that is, if a sequence $\{u_n\}\subset W_{1}$ satisfying
 $$u_n\rightharpoonup u \quad \mbox{ and } \quad \limsup_{n\rightarrow\infty}\langle\mathcal{F}(u_n),u_n-u\rangle\leq 0,$$
 then $u_n\rightarrow u$ in $W_{1}$.
 Similarly, we can also obtain that $v_n\rightarrow v$ in $W_{2}$. Therefore,
 $\{(u_n,v_n)\}\rightarrow (u,v)$ in $W$. \qed
 \vskip2mm
 \noindent
{\bf Proof of Theorem 3.1.} It is obvious that $I(0)=0$.  By Lemma 3.5, Lemma 3.6 and Lemma 3.8, all conditions of Lemma 2.9 hold. Then system \eqref{eq1} possesses a nontrivial weak solution which is a critical point of $I$. \qed

\subsection{Multiplicity}
In this section, by using the Symmetric Mountain Pass Theorem, we can obtain the following multiplicity result.\\
 \noindent
{\bf Theorem 3.9. }
  Assume that $(\phi_1)$--$(\phi_4)$,  $(V_0)$, $(V_1)$, $(M_0)$, $(M_1)$, $(F_0)$, $(F_1)$,  $(F_4)$ and the following conditions hold:
 \begin{itemize}
 \item[ $(F_{3}^{\prime})$]
 $$
 \lim_{|(u,v)|\rightarrow+\infty}\frac{F(x,u,v)}{|u|^{m_1}+|v|^{m_2}}=+\infty \mbox{\ \ uniformly in\ \ } x \in \mathbb{R}^N;
 $$
 \item[$(F_6)$] $F(x,-u,-v)=F(x,u,v)$ \; for all $(x,u,v)\in \mathbb{R}^N\times \mathbb{R}\times \mathbb{R}.$
  \end{itemize}
  Then system \eqref{eq1} possesses infinitely many weak solutions $\{(u_k,v_k)\}$ such that
 \begin{eqnarray*}
        I(u_k,v_k):
= &       \widehat{M_{1}}\left(\int_{\mathbb{R}^N}\Phi_1(|\nabla u_k|)dx\right)
          + \widehat{M_{2}}\left(\int_{\mathbb{R}^N}\Phi_2(|\nabla v_k|)dx\right)
          + \int_{\mathbb{R}^N}V_{1}(x)\Phi_1(| u_k|)dx
          \\
&~~~
          + \int_{\mathbb{R}^N}V_{2}(x)\Phi_2(|v_k|)dx
          - \int_{\mathbb{R}^{N}}F(x,u_k,v_k)dx
\rightarrow + \infty \;\; \mbox{ as } k \rightarrow \infty.
 \end{eqnarray*}

\vskip0mm
\par
To apply the Symmetric Mountain Pass Theorem (i.e., Lemma 2.10), we need the following knowledge. One can see the details in  \cite{wang2017,J.F. Zhao1991,N.T. Chung2013}.
 Since $W$ is a reflexive and separable Banach spaces,  there exist two sequences  $\{e_{ij}: j\in \mathbb{N}\}\subset {W_{i}}~(i=1,2)$ and $\{e_{ij}^{*}: j\in {\mathbb{N}\}\subset {W_{i}}}^{*}~(i=1,2)$ such that
 \begin{equation*}\label{3.2.1}
 W_{i}=\overline{\mbox {span}\{e_{ij}: j=1,2,\cdots\}}, \quad {W_{i}}^{*}=\overline{\mbox{span}\{e_{ij}^{*}: j=1,2,\cdots\}},\quad i=1,2,
 \end{equation*}
 and
 \begin{equation*}\label{3.2.2}
 e^{*}_{in}(e_{im})=
 \begin{cases}
  \begin{array}{ll}
 1 &\mbox{if } n=m,\\
 0 &\mbox{if } n\neq m,\\
    \end{array}
 \end{cases}i=1,2.
 \end{equation*}
 Let $Y_{i(k)}$ and $Z_{i(k)}$ be the subsets of $W_{i}$ defined by
 \begin{equation*}\label{3.2.3}
 Y_{i(k)}:=\mbox{span}\{e_{ij}: j=1,\cdots,k\}, \quad Z_{i(k)}:=\overline{\mbox{span}\{e_{ij}: j=k+1,\cdots\}}, \quad i=1,2.
 \end{equation*}
  Then
 $$W_{i}=Y_{i(k)}\oplus Z_{i(k)},\quad i=1,2, \quad k\in\mathbb{N}.$$

 Moreover, since the embeddings  $W_{i}\hookrightarrow L^{\Psi_i}(\mathbb{R}^{N})~(i=1,2)$ and $W_{i}\hookrightarrow L^{l_i}(\mathbb{R}^{N})~(i=1,2)$ are compact,  with a similar discussion as \cite{wang2017,N.T. Chung2013}, we can get
 \begin{equation}\label{3.2.4}
 \alpha_{i(k)}:=\mbox{sup} \left\{\|z\|_{\Psi_{i}}: \|z\|_{i}=1, z\in Z_{i(k)}\right\}\rightarrow 0
 \end{equation}
and
\begin{equation}\label{3.2.5-}
 \beta_{i(k)}:=\mbox{sup} \left\{\|z\|_{L^{l_i}}: \|z\|_{i}=1, z\in Z_{i(k)}\right\}\rightarrow 0, \quad i=1, 2, \mbox{ as } k\rightarrow \infty.
 \end{equation}

In addition, for Banach space $W=W_{1}\times W_{2}$, there exists a sequence $\{\eta_{(j)}\}\subset W$ defined by
 \begin{equation*}\label{3.2.6}
 \eta_{(j)}=
 \begin{cases}
  \begin{array}{ll}
 (e_{1n},0) &\mbox{if }j=2n-1,\\
 (0,e_{2n}) &\mbox{if }j=2n,  \quad      \mbox{ for } n\in \mathbb{N},\\
    \end{array}
 \end{cases}
 \end{equation*}
 such that
 \begin{itemize}
 	\item[$\rm(1)$]
 $$W=\overline{\mbox{span}\{\eta_{(j)}: j=1,2,\cdots\}},$$
 \item[$\rm(2)$]
 $$W=Y_k\oplus Z_k,$$
 where
 $$Y_k:=\text{span}\{\eta_{(j)}: j=1,\cdots,k\} \quad \text{and} \quad Z_k:=\overline{\mbox{span}\{\eta_{(j)}: j=k+1,\cdots\}}.$$
 \end{itemize}

\vskip2mm
\noindent
 {\bf Lemma 3.10. }
  Suppose that $(\phi_1)$--$(\phi_3)$, $(M_0)$, $(V_0)$ and $(F_1)$ hold. Then there are two constants $\rho, \alpha>0$ and $k\in\mathbb{N}$ such that
  $I\mid _{\partial B_{\rho}\cap Z_{2k}}\geq \alpha$.

  \vskip2mm
  \noindent
  {\bf Proof. }
 For $(u,v)\in Z_{2k}$,
in view of the  inequality \eqref{3.1.2-1}, \eqref{3.1.6}, \eqref{3.1.7}, \eqref{3.2.4}, \eqref{3.2.5-}, $(M_{0})$, Young's inequality,  Lemma 2.4,  Remark 3.2  and the inequality (66) in \cite{Xie2018}, we have
 \begin{align*}
          I(u,v)
 &  =     \widehat{M}_1\left(\int_{\mathbb{R}^N}\Phi_1(|\nabla u|)dx\right)
          +  \widehat{M}_2\left(\int_{\mathbb{R}^N}\Phi_2(|\nabla v|)dx\right)
          +  \int_{\mathbb{R}^N}\left(    V_{1}(x)\Phi_1(| u|)
                                         +  V_{2}(x)\Phi_2(| v|)
                                         -  F(x,u,v)
                                  \right)
          dx
          \\
 &\geq       C_{1,3}\int_{\mathbb{R}^N}\Phi_1(|\nabla u|)dx
          +  C_{2,3}\int_{\mathbb{R}^N}\Phi_2(|\nabla v|)dx
          + \int_{\mathbb{R}^N}V_{1}(x)\Phi_1(| u|)dx
          +  \int_{\mathbb{R}^N}V_{2}(x)\Phi_2(| v|)dx
          \nonumber\\
&~~~~
          - C_{5} \int_{\mathbb{R}^N}|u|^{l_{1}}dx
          - C_{5} \int_{\mathbb{R}^N}|v|^{l_{2}}dx
          - C_{5} \int_{\mathbb{R}^N}\Psi_1(u)dx
          - C_{5} \int_{\mathbb{R}^N}\Psi_2(v)dx
                  \\
&\geq       C_{1,3}
            \min\left\{\|\nabla u\|_{\Phi_1}^{l_1},\|\nabla u\|_{\Phi_1}^{m_1}\right\}
          +
            \min\left\{\|u\|_{\Phi_1,V_{1}}^{l_1},\|u\|_{\Phi_1,V_{1}}^{m_1}\right\}
          - C_{5}\|u\|_{L^{l_{1}}}^{l_{1}}
          - C_{5}\max\left\{\|u\|_{\Psi_1}^{l_{\Psi_1}},\|u\|_{\Psi_1}^{m_{\Psi_1}}\right\}
          \nonumber\\
&~~
          + C_{2,3}
            \min\left\{\|\nabla v\|_{\Phi_2}^{l_2},\|\nabla v\|_{\Phi_2}^{m_2}\right\}
          +
            \min\left\{\|v\|_{\Phi_2,V_{2}}^{l_2},\|v\|_{\Phi_2,V_{2}}^{m_2}\right\}
         - C_{5}\|v\|_{L^{l_{2}}}^{l_{2}}
         - C_{5}\max\left\{\|v\|_{\Psi_2}^{l_{\Psi_2}},\|v\|_{\Psi_2}^{m_{\Psi_2}}\right\}
          \nonumber\\
&\geq
       \frac{\min\{C_{1,3},1\}}{2^{m_{1}-1}}
       \|u\|_{1}^{l_1}
       +
        \frac{\min\{C_{2,3},1\}}{2^{m_{2}-1}}
       \|v\|_{2}^{l_2}
        -\min\{C_{1,3},1\}
        -\min\{C_{2,3},1\}
          - C_{5}\beta_{1(k)}^{l_{1}}\|u\|_{1}^{l_{1}}
             \nonumber\\
&~~~~
          - C_{5}\beta_{2(k)}^{l_{2}}\|v\|_{2}^{l_{2}}
          - C_{5}\max\left\{\alpha_{1(k)}^{l_{\Psi_1}}\|u\|_{1}^{l_{\Psi_1}},
                            \alpha_{1(k)}^{m_{\Psi_1}}\|u\|_{1}^{m_{\Psi_1}}
                    \right\}
          - C_{5}\max\left\{\alpha_{2(k)}^{l_{\Psi_2}}\|v\|_{2}^{l_{\Psi_2}},
                            \alpha_{2(k)}^{m_{\Psi_2}}\|v\|_{2}^{m_{\Psi_2}}
                     \right\}.
 \end{align*}
 Since $\alpha_{i(k)}\rightarrow0,~\beta_{i(k)}\rightarrow0,~i=1,2,$ as $k\rightarrow\infty$, then above inequality implies that for  a large constant $\rho>0$, there exists a large $k\in \mathbb{N}$ such that
 ${I}\mid _{\partial B_{\rho}\cap Z_{2k}}\geq \alpha$ for some $\alpha>0$. \qed
 \vskip2mm
\noindent
 {\bf Lemma 3.11. }
  Suppose that $(\phi_1)$--$(\phi_3)$,  $(M_1)$, $(V_0)$, $(F_1)$ and $(F_{3}^{\prime})$ hold. Then for each finite dimensional subspace $\widetilde{W}\subset W$, there exists a positive constant $R=R(\widetilde{W})$ such that $I\leq 0$ on $\widetilde{W}\backslash B_{R(\widetilde{W})}$.

\vskip2mm
\noindent
  {\bf Proof. } For each finite dimensional subspace $\widetilde{W}\subset W$,
 one has $\widetilde{W}\subseteq W_{1,1}\times W_{2,1},$ where $W_{1,1}$ and $W_{2,1}$ are finite dimensional subspaces of $W_{1}$ and $W_{2}$, respectively.
 Since any two norms in finite dimensional space are equivalent,  there exist positive constants $h_1,\;h_2,\;h_3,\;h_4,\;h_5,\;h_6,\;h_7$ and $h_8$  such that
 \begin{align}
 &h_1\|u\|_{1}\leq\|u\|_{L^{m_1}}\leq h_2\|u\|_{1},\quad \forall  u\in W_{1,1},\label{3.2.5}\\
 &h_3\|v\|_{2}\leq\|v\|_{L^{m_2}}\leq h_4\|v\|_{2},\quad \forall  v\in W_{2,1},\label{3.2.6}\\
  &h_5\|u\|_{1}\leq\|u\|_{L^{l_1}}\leq h_6\|u\|_{1},\quad \forall  u\in W_{1,1},\label{3.2.5-}\\
 &h_7\|v\|_{2}\leq\|v\|_{L^{l_2}}\leq h_8\|v\|_{2},\quad \forall  v\in W_{2,1},\label{3.2.6-}
 \end{align}
 where $h_1, h_2, h_5, h_6$ and $h_3, h_4, h_7, h_8$ depend on the spatial dimension of $W_{1,1}$ and $W_{2,1}$, respectively.
 Moreover, $(F_1)$, $(F_{3}^{\prime})$ and the continuity of $F$ imply that for any given constant
 $$
 L> \max\bigg\{\frac{4\max\{C_{1,4},1\}+C(L)h_{6}^{l_1}}{h_{1}^{m_1}},
            \frac{4\max\{C_{2,4},1\}+C(L)h_{8}^{l_2}}{h_{3}^{m_2}}
     \bigg\},
     $$ there exists a constant $C(L)>0$ such that
\begin{equation}\label{3.2.--}
 F(x,u,v)\geq L(|u|^{m_1}+|v|^{m_2})-C(L)(|u|^{l_1}+|v|^{l_2}), \quad \forall  (x,u,v)\in \mathbb{R}^{N}\times \mathbb{R}\times \mathbb{R}.
 \end{equation}
  Then,  by
 $(M_0)$,  \eqref{3.1.6}, \eqref{3.1.7},   \eqref{3.2.5}-\eqref{3.2.--}, Lemma 2.4 and inequality (67) in \cite{Xie2018},  we have
 \begin{align*}
            I(u,v)
&   =      \widehat{M}_{1}\left(\int_{\mathbb{R}^N}\Phi_1(|\nabla u|)dx\right)
          + \widehat{M}_{2}\left(\int_{\mathbb{R}^N}\Phi_2(|\nabla v|)dx\right)
          + \int_{\mathbb{R}^N}V_{1}(x)\Phi_1(|u|)dx
           \\
&~~~
          + \int_{\mathbb{R}^N}V_{2}(x)\Phi_2(|v|)dx
          - \int_{\mathbb{R}^N}F(x,u,v)dx
          \\
&\leq
            \left(  C_{1,4}\int_{\mathbb{R}^N}\Phi_1(|\nabla u|)dx
                  +  \int_{\mathbb{R}^N}V_{1}(x)\Phi_1(|u|)dx
           \right)
          -  L\int_{\mathbb{R}^N}|u|^{m_1}dx
          +  C(L)\int_{\mathbb{R}^N}|u|^{l_1}dx
          \\
&~~~
          +
           \left(
                    C_{2,4}\int_{\mathbb{R}^N}\Phi_2(|\nabla v|)dx
                  + \int_{\mathbb{R}^N}V_{2}(x)\Phi_2(|v|)dx
           \right)
          -  L\int_{\mathbb{R}^N}|v|^{m_2}dx
          +  C(L)\int_{\mathbb{R}^N}|v|^{l_2}dx
          \\
 &\leq              C_{1,4}\|\nabla u\|_{\Phi_{1}}^{l_1}+C_{1,4}\|\nabla u\|_{\Phi_{1}}^{m_1}
                 +\|u\|_{\Phi_{1},V_{1}}^{l_1}+\|u\|_{\Phi_{1},V_{1}}^{m_1}
                 -  L\|u\|_{L^{m_1}}^{m_1}
                 +  C(L)\|u\|_{L^{l_1}}^{l_1}
                 \\
&~~~
     +   C_{2,4}\|\nabla v\|_{\Phi_{2}}^{l_2}+C_{2,4}\|\nabla v\|_{\Phi_{2}}^{m_2}
                  +\|v\|_{\Phi_{2},V_{2}}^{l_2}+\|v\|_{\Phi_{2},V_{2}}^{m_2}
                 -  L\|v\|_{L^{m_2}}^{m_2}
                 +  C(L)\|v\|_{L^{l_2}}^{l_2}
          \\
 &\leq           2\max\{C_{1,4},1\}\|u\|_{1}^{l_1}
                 +2\max\{C_{1,4},1\}\|u\|_{1}^{m_1}
                 +4\max\{C_{1,4},1\}
                 -  Lh_{1}^{m_1}\|u\|_{1}^{m_1}
                 + C(L)h_{6}^{l_1}\|u\|_{1}^{l_1}
                 \\
&~~~
     +             2\max\{C_{2,4},1\}\|v\|_{2}^{l_2}
     +             2\max\{C_{2,4},1\}\|v\|_{2}^{m_2}
                  +4\max\{C_{2,4},1\}
                 -  Lh_{3}^{m_2}\|v\|_{2}^{m_2}
                 + C(L)h_{8}^{l_2}\|v\|_{2}^{l_2}
          \\
 &=              \left(2\max\{C_{1,4},1\}+C(L)h_{6}^{l_1}\right)\|u\|_{1}^{l_1}
                 -  \left(Lh_{1}^{m_1}-2\max\{C_{1,4},1\}\right)\|u\|_{1}^{m_1}
                 +4\max\{C_{1,4},1\}
                 \\
&~~~
     +          \left(2\max\{C_{2,4},1\}+C(L)h_{8}^{l_2}\right)\|v\|_{2}^{l_2}
                 -  \left(Lh_{3}^{m_2}-2\max\{C_{2,4},1\}\right)\|v\|_{2}^{m_2}
                  +4\max\{C_{2,4},1\}.
\end{align*}
Note that $l_i\leq m_i$ $(i=1, 2)$. Then the above inequality implies that
$$\lim_{r\rightarrow \infty}\sup_{(u,v)\in \partial B_r\cap\widetilde{W}}I(u,v)=-\infty.$$
Thus, there exists an $R=R(\widetilde{W})$ such that $I\leq 0$ on $\widetilde{W}\backslash B_{R(\widetilde{W})}$. \qed

 \vskip2mm
 \noindent
{\bf Proof of Theorem 3.9.} By $(F_6)$, it is obvious that $I$ is even in $W$.  By Lemma 3.10, Lemma 3.11 and Lemma 3.8, all conditions of Lemma 2.10 hold. Then system \eqref{eq1} possesses infinitely many weak solutions $\{(u_k,v_k)\}$ and  $I(u_k,v_k)\to +\infty$ as $k\to+\infty$. \qed

\section{Results for the scalar equation}
In this section, we study the existence and multiplicity of solutions for the following generalized Kirchhoff elliptic equation in Orlicz-Sobolev spaces:
 \begin{equation}\label{eq5}
 \left\{
  \begin{array}{ll}
 -M\left(\int_{\mathbb{R}^{N}}\Phi(|\nabla u|)dx\right)\Delta_{\Phi}u+V(x)\phi(|u|)u=f(x,u), &x\in \mathbb{R}^N,\\
  u\in W^{1,\Phi_{1}}(\mathbb{R}^N),
    \end{array}
 \right.
 \end{equation}
where  $\phi:(0,+\infty)\rightarrow(0,+\infty)$ is a function
which satisfies:
\begin{itemize}
	 \item[$(\phi_1)'$] $\phi\in C^1(0,+\infty)$, $t\phi(t)\rightarrow0$ as
 $t\rightarrow0$, $t\phi(t)\rightarrow+\infty$ as
 $t\rightarrow+\infty$;
 \item[$(\phi_2)'$] $t\rightarrow t\phi(t)$ are strictly increasing;
 \item[$(\phi_3)'$] $1<l:=\inf_{t>0}\frac{t^2\phi(t)}{\Phi(t)}\leq\sup_{t>0}\frac{t^2\phi(t)}{\Phi(t)}=:m<\min\{N,l^{\ast}\}$, where $\Phi(t):=\int_{0}^{|t|}s\phi(s)ds, \ t\in \mathbb{R}$, $l^{\ast}=\frac{lN}{N-l}$;
 \item[$(\phi_4)'$]
    there exist positive constants $c_{1}$ and $c_{2}$  such that
    $$
    c_{1}|t|^{l}\le \Phi(t)\le c_{2}|t|^{l},\ \ \forall |t|<1;
    $$
\end{itemize}
 Moreover, we introduce the following conditions on $f$, $V$ and $M$:
\begin{itemize}
\item[$(F_0)'$]  $f: \mathbb{R}^{N}\times \mathbb{R}\rightarrow \mathbb{R}$ is a $C^1$ function such that
                $f(x,0)=0, \ x\in \mathbb{R}^{N}$;
\end{itemize}
\begin{itemize}
\item[$(V_0)'$]  $V \in C(\mathbb{R}^{N},\mathbb{R})$ and $\inf_{\mathbb{R}^{N}}V(x)> 1$;
\end{itemize}
\begin{itemize}
\item[$(V_1)'$]  there exist  constants $c_{3}>0$ such that
                $$
                   \lim_{|z|\rightarrow \infty} \mbox{meas}\{x \in \mathbb{R}^{N}: |x-z|\leq c_{3},
                          V(x)\leq c_{4}\}=0\;\;\mbox{for every}\; c_{4}>0,
                $$
                where meas$(\cdot)$ denotes the Lebesgue measure in $\mathbb{R}^{N}$;
\end{itemize}
\begin{itemize}
\item[$(M_0)'$]  $M \in C(\mathbb{R}^+,\mathbb{R}^+) $  and $c_{5} \leq M(t)\leq c_{6},\; \forall \;t\geq 0$ for some $c_{5}, c_{6}>0$;
\item[$(M_1)'$]  $\widehat{M}(t):=\int_{0}^{t}M(s)ds\geq M(t)t$.
\end{itemize}
\par
Similar to the results in  section 3, we can get the following results.
\vskip2mm
\noindent
{\bf Theorem 4.1. }
  Assume that $\phi$ and $f$ satisfy $(\phi_1)'$-$(\phi_4)'$, $(F_0)'$, $(M_0)'$, $(M_1)'$, $(V_0)'$, $(V_1)'$  and the following conditions:
\begin{itemize}
 \item[$(F_1)'$] there exist a constant $c_{7}>0$ and a continuous function $\psi: [0,+\infty)\rightarrow \mathbb{R}$ such that
 $$|f(x, u)|\leq c_{7}(|u|^{l-1}+\psi(|u|))$$
 for all $(x,u)\in \mathbb{R}^{N}\times\mathbb{R}$,
 where
 $$\Psi(t):=\int_{0}^{|t|}\psi(s)ds, \quad t\in \mathbb{R}$$
 is an $N$-function  satisfying
 $$m<l_{\Psi}:=\inf_{t>0}\frac{t\psi(t)}{\Psi(t)}\leq\sup_{t>0}\frac{t\psi(t)}{\Psi(t)}
 =:m_{\Psi}<l^{*}:=\frac{lN}{N-l};$$
 \item[$(F_2)'$] there exists a constant $c_{8}\in [0,1)$ such that
 $$\limsup_{u\rightarrow 0}\frac{F(x, u)}{\Phi(u)}=c_{8}\mbox{\ \ uniformly in\ \ } x \in \mathbb{R}^{N},$$
 where $F(x,u)=\int_0^uf(x,s)ds$ for all $(x,u)\in \mathbb R^N\times\mathbb R;$
 \item[$(F_3)'$] there exists a domain  $G\subset \mathbb{R}^{N}$ such that
 $$\lim_{u\rightarrow \infty}\frac{F(x, u)}{|u|^{m}}=+\infty,\mbox{\ \ for a.e.\ \ } x \in G;$$
 \item[$(F_4)'$] there exists a continuous function $\overline{\gamma}: [0,\infty)\rightarrow \mathbb{R}$ and constants
$  \sigma
\in
   \bigg[
                     \frac{l(m_{\overline{\Gamma}}-1)}{m_{\overline{\Gamma}}},
          \min\bigg\{l,
                    \frac{l^{\ast}(l_{\overline{\Gamma}}-1)}{l_{\overline{\Gamma}}}
              \bigg\}
   \bigg),
$
$c_9, r_{2}>0$
 such that
 \begin{equation}\label{5.1.3}
      \overline{\Gamma}\left(\frac{F(x,u)}{|u|^{\sigma_1}}\right)
 \leq c_9\overline{F}(x,u), \quad \mbox{ for all } x\in \mathbb{R}^N \mbox{ and all } u\in \mathbb{R} \mbox{ with }|u|\geq r_{2},
 \end{equation}
 where  $\overline{\Gamma}(t):=\int_{0}^{|t|}\overline{\gamma}(s)ds,~ t\in\mathbb{R}$, is an $N$-function with
 \begin{equation}\label{5.1.4}
       1
  <    l_{\overline{\Gamma}}:=\inf_{t>0}\frac{t\overline{\gamma}(t)}{\overline{\Gamma}(t)}
  \leq \sup_{t>0}\frac{t\overline{\gamma}(t)}{\overline{\Gamma}(t)}=:m_{\overline{\Gamma}}
  <+\infty
 \end{equation}
 and
 $$\overline{F}(x,u):=\frac{1}{m}f(x,u)u-F(x,u), ~~~\forall  (x,u)\in \mathbb{R}^N\times \mathbb{R}.$$
\end{itemize}
 Then the  equation \eqref{eq5} has a nontrivial weak solution.

 \vskip2mm
 \noindent
 {\bf Theorem 4.2. }
  Assume that $(\phi_1)'$-$(\phi_4)'$, $(F_0)'$, $(F_1)'$, $(M_0)'$, $(M_1)'$, $(V_0)'$, $(V_1)'$,      $(F_4)'$  and the following conditions hold:
 \begin{itemize}
 \item[ $(F_{3}^{\prime})'$]
 $$
 \lim_{t\rightarrow+\infty}\frac{F(x,u)}{|u|^{m}}=+\infty \mbox{\ \ uniformly in\ \ } x \in \mathbb{R}^N;
 $$
\item[ $(F_5)'$] $F(x,-u)=F(x,u)$, for all $(x,u)\in \mathbb{R}^N\times \mathbb{R}.$
  \end{itemize}
  Then the equation \eqref{eq5} possesses infinitely many weak solutions $\{u_k\}$ such that
  \begin{eqnarray*}
        I(u_k):
=       \widehat{M}\left(\int_{\mathbb{R}^N}\Phi(|\nabla u_k|)dx\right)
          + \int_{\mathbb{R}^N}V(x)\Phi(| u_k|)dx
          - \int_{\mathbb{R}^{N}}F(x,u_k)dx
\rightarrow + \infty ,\;\; \mbox{ as } k \rightarrow \infty.
 \end{eqnarray*}

\vskip2mm
\noindent
{\bf Remark 4.3. }
\begin{itemize}
\item[(\uppercase\expandafter{\romannumeral1})]
It is clear that the conditions $(F_3)'$ and $(F_4)'$ extend the conditions $(F_3)$ and $(F_5)$  in \cite{Tang2019};
\item[(\uppercase\expandafter{\romannumeral2})]
Comparing Theorem 1.1 in \cite{Liu-Shibo2019} and Theorem 1.5 in \cite{Alves2014} with  Theorem 4.1, it is easy to see that the condition $(F_3)'$ is weaker than the following global (A-R) condition:\\
(A-R) there exist $\theta>m$ such that for all $u \in \mathbb R/\{0\}$,
$$
0<F(u):=\int_{0}^{u}f(x,s)ds \leq \frac{1}{\theta}uf(x,u).
$$
\item[(\uppercase\expandafter{\romannumeral3})] If we consider the system (\ref{eq5}) on a bounded domain $\Omega$ with Dirichlet boundary condition, then it is natural that we restrict those assumptions of Theorem 4.1 on the bounded domain $\Omega$.
Then the condition $(F_4)'$ is different from the
condition $(f_4)$ in Theorem 5.1 of \cite{wang2017} and $(f_2)$ in  \cite{Carvalho2015}. To exemplify this, let $F(x,t)= \left(\sin(2\pi x_{1})+|\sin(2\pi x_{1})|\right)\left(  |t|^{5}\ln(1+|t|)\right)$. If we choose  $l_{\overline{\Gamma}}= \frac{6}{5}$, then the condition $(F_4)$ holds. But in \cite{Carvalho2015,wang2017}, the constant $l_{\Gamma}$ must satisfy $l_{\Gamma}>\frac{N}{l_{i}}$, i.e., $l_{\Gamma}>\frac{3}{2}>\frac{6}{5}$. So the condition $(f_4)$ in \cite{wang2017} and $(f_2)$ in  \cite{Carvalho2015} do not hold.
\end{itemize}

\section{Example}\label{section 6}
\begin{eqnarray}\label{eq6}
 \begin{cases}
  -M_{1}\left(\int_{\mathbb{R}^{6}}\left(|\nabla u|^{4}+|\nabla u|^{5}\right)dx\right)\mbox{div}[(4|\nabla u|^{2}+5|\nabla u|^{3})\nabla u]+V_{1}(x)(4|u|^{2}+5|u|^{3})u\\
 \;\;\; = F_u(x,u,v), \ \ x\in \mathbb{R}^6,\\
   -M_{2}\left(\int_{\mathbb{R}^{6}}\left(|\nabla v|^{4}\ln(e+|\nabla v|)\right)dx\right)
   \mbox{div}\left[\left(4|\nabla v|^{2}\ln(e+|\nabla v|)+\dfrac{|\nabla v|^{3}}{e+|\nabla v|}\right)\nabla v\right]\\
 \;\;\; +V_{1}(x)\left(4|v|^{2}\ln(e+| v|)+\dfrac{|v|^{3}}{e+|v|}\right)v
= F_v(x,u,v), \ \ x\in \mathbb{R}^6,\\
 \end{cases}
\end{eqnarray}
where
\begin{eqnarray}
&M_{1}(t)=2+\frac{1}{e+t},\;t\geq 0,\; M_{2}(t)=3+\frac{1}{et^{\frac{2}{3}}+3t},\;t\geq 0, \label{6.1.1}\\
&F(x,t,s)= \left(\sin(2\pi x_1)+|\sin(2\pi x_1)|\right)
           \left(  |t|^{5}\ln(1+|t|)+|s|^{5}\ln(1+|s|)+|t|^{3}|s|^{3}\right).\label{6.1.2}
\end{eqnarray}
Let $N=6$, $\phi_1(t)=4|t|^{2}+5|t|^{3}$ and $\phi_2(t)=4|t|^{2}\ln(e+| t|)+\dfrac{|t|^{3}}{e+|t|}$. Then $\phi_i (i=1,2)$ satisfy $(\phi_1)$-$(\phi_4)$, $M_{i}$ satisfy $(M_0)$-$(M_1)$ by (\ref{6.1.1}), $l_{1}=l_{2}=4$, $m_{1}=m_{2}=5$ and $\Phi_1(t)=|t|^{4}+|t|^{5}$, $\Phi_2(t)=|t|^{4}\ln(e+|t|)$. So, $l^{*}_{1}=l^{*}_{2}=12$.
\par
 Let $V_1(x)=\sum_{i=1}^{6}x_i^{2}+1$ and $V_2(x)=\sum_{i=1}^{6}x_i^{3}+2$
 for all $(x,t,s)\in \mathbb{R}^{N}\times \mathbb{R}\times \mathbb{R}$. Then it is obvious that $V_i, i=1,2$ satisfy (V0) and (V1).
 \par
By (\ref{6.1.2}), we have
\begin{eqnarray}
&F_{t}(x,t,s)= \left(\sin(2\pi x_1)+|\sin(2\pi x_1)|\right)
            \left(    5|t|^{3}t\ln(1+|t|)
              +\frac{|t|^{4}t}{1+|t|}
              + 3|t||s|^{3}t
              \right), \label{6.1.3}\\
&F_{s}(x,t,s)= \left(\sin(2\pi x_1)+|\sin(2\pi x_1)|\right)
           \left(     5|s|^{3}s\ln(1+|s|)
              +\frac{|s|^{4}s}{1+|s|}
              + 3|t|^{3}|s|s
              \right).\label{6.1.4}
\end{eqnarray}
Hence
\begin{align}\label{6.1.5}
\overline{F}(x,t,s)
=&     \left(\sin(2\pi x_1)+|\sin(2\pi x_1)|\right)
      \left(    \frac{|t|^{6}}{5(1+|t|)}
              + \frac{|s|^{6}}{5(1+|s|)}
              + \frac{1}{5}|t|^{3}|s|^{3}
      \right)\nonumber\\
\geq &  \left(\sin(2\pi x_1)+|\sin(2\pi x_1)|\right)
        \left( \frac{|t|^{6}}{5(1+|t|)}+\frac{|s|^{6}}{5(1+|s|)}\right)\nonumber\\
\geq &
\begin{cases}
             \frac{\sin(2\pi x_1)+|\sin(2\pi x_1)|}{10}
        (|t|^{5}+|s|^{5}),
                   \;\;\; \;\;\; \;\; \;             \text { if }\;\;|t|\geq 1, |s|\geq 1,\\
             \frac{\sin(2\pi x_1)+|\sin(2\pi x_1)|}{10}
        (|t|^{6}+|s|^{5}),
                   \;\;\; \;\;\; \;\; \;             \text { if }\;\;0\leq |t|< 1, |s|\geq 1,\\
             \frac{\sin(2\pi x_1)+|\sin(2\pi x_1)|}{10}
        (|t|^{5}+|s|^{6}),
                   \;\;\; \;\;\; \;\; \;             \text { if }\;\;|t|\geq 1, 0\leq |s|< 1,\\
              \frac{\sin(2\pi x_1)+|\sin(2\pi x_1)|}{10}
        (|t|^{6}+|s|^{6}),
                   \;\;\; \;\;\; \;\; \;             \text { if }\;\;0\leq |t|< 1, 0\leq |s|< 1.
      \end{cases}
\end{align}
It is easy to see that conditions $(F_0)$ and $(F_5)$ hold.
Since
$$
\lim_{|(t,s)|\rightarrow 0}\frac{F(x,t,s)}{|t|^{m_{1}}+|s|^{m_{2}}}=0, \;\;\mbox{and }
\lim_{|(t,s)|\rightarrow +\infty}\frac{F(x,t,s)}{|t|^{m_{1}}+|s|^{m_{2}}}=+\infty,
$$
by (2) of Lemma 2.3,
we can see that $(F_2)$  and $(F_3)$ hold with $G=(1/8,3/8)\times \mathbb{R}^{5}$.
Choose  $\Psi_1(t)=\Psi_2(t)=t^{6}, \overline{\Gamma}(t)=|t|^{\frac{6}{5}}$ and $\sigma_{1}, \sigma_{2}=\frac{11}{6}$. Then
{\small \begin{align*}
&   \limsup_{|(u,v)|\rightarrow \infty}
      \left(\frac{|F(x,u,v)|}{|u|^{\frac{11}{6}}+|v|^{\frac{11}{6}}}\right)^{\frac{6}{5}}\frac{1}{\overline{F}(x,u,v)}\\
&
\leq
\begin{cases}
             \limsup_{|(u,v)|\rightarrow \infty}
 \frac{  10
       \left(\sin(2\pi x_1)+|\sin(2\pi x_1)|\right)^{\frac{11}{5}}
       \left(|u|^{5}\ln(1+|u|)+|v|^{5}\ln(1+|v|)+|u|^{3}|v|^{3}\right)^{\frac{6}{5}}
      }
        {\left(|u|^{\frac{11}{6}}+|v|^{\frac{11}{6}}\right)^{\frac{6}{5}}(|u|^{5}+|v|^{5})},
                   \;\;\; \;\;\; \;\; \;             \text { if }\;\;|u|\geq 1, |v|\geq 1,\\
            \limsup_{|(u,v)|\rightarrow \infty}
 \frac{  10
       \left(\sin(2\pi x_1)+|\sin(2\pi x_1)|\right)^{\frac{11}{5}}
       \left(|u|^{5}\ln(1+|u|)+|v|^{5}\ln(1+|v|)+|u|^{3}|v|^{3}\right)^{\frac{6}{5}}
      }
        {\left(|u|^{\frac{11}{6}}+|v|^{\frac{11}{6}}\right)^{\frac{6}{5}}(|u|^{6}+|v|^{5})},
                   \;\;\; \;\;\; \;\; \;             \text { if }\;\;0\leq |u|< 1, |v|\geq 1,\\
             \limsup_{|(u,v)|\rightarrow \infty}
 \frac{  10
       \left(\sin(2\pi x_1)+|\sin(2\pi x_1)|\right)^{\frac{11}{5}}
       \left(|u|^{5}\ln(1+|u|)+|v|^{5}\ln(1+|v|)+|u|^{3}|v|^{3}\right)^{\frac{6}{5}}
      }
        {\left(|u|^{\frac{11}{6}}+|v|^{\frac{11}{6}}\right)^{\frac{6}{5}}(|u|^{5}+|v|^{6})},
                   \;\;\; \;\;\; \;\; \;             \text { if }\;\;|u|\geq 1, 0\leq |v|< 1
      \end{cases} \\
&
\leq
\begin{cases}
             40\left(\sin(2\pi x_1)+|\sin(2\pi x_1)|\right)^{\frac{11}{5}}
   \limsup_{|(u,v)|\rightarrow \infty}
   \frac{|u|^{6}(\ln(1+|u|))^{\frac{6}{5}}+|v|^{6}(\ln(1+|v|))^{\frac{6}{5}}+|u|^{\frac{36}{5}}+|v|^{\frac{36}{5}}}
        {|u|^{\frac{36}{5}}+|v|^{\frac{36}{5}}},
                   \;\;\;              \text { if }\;\;|u|\geq 1, |v|\geq 1,\\
            40\left(\sin(2\pi x_1)+|\sin(2\pi x_1)|\right)^{\frac{11}{5}}
   \limsup_{|(u,v)|\rightarrow \infty}
   \frac{|u|^{6}(\ln(1+|u|))^{\frac{6}{5}}+|v|^{6}(\ln(1+|v|))^{\frac{6}{5}}+|u|^{\frac{36}{5}}+|v|^{\frac{36}{5}}}
        {|u|^{\frac{41}{5}}+|v|^{\frac{36}{5}}},
                   \;\;\;           \text { if }\;\;0\leq |u|< 1, |v|\geq 1,\\
            40\left(\sin(2\pi x_1)+|\sin(2\pi x_1)|\right)^{\frac{11}{5}}
   \limsup_{|(u,v)|\rightarrow \infty}
   \frac{|u|^{6}(\ln(1+|u|))^{\frac{6}{5}}+|v|^{6}(\ln(1+|v|))^{\frac{6}{5}}+|u|^{\frac{36}{5}}+|v|^{\frac{36}{5}}}
        {|u|^{\frac{36}{5}}+|v|^{\frac{41}{5}}},
                   \;\;\;             \text { if }\;\;|u|\geq 1, 0\leq |v|< 1
      \end{cases}
\\
&
\leq
\begin{cases}
             40\left(\sin(2\pi x_1)+|\sin(2\pi x_1)|\right)^{\frac{11}{5}}
   \limsup_{|(u,v)|\rightarrow \infty}
   \frac{|u|^{6}|u|^{\frac{6}{5}}+|v|^{6}|v|^{\frac{6}{5}}+|u|^{\frac{36}{5}}+|v|^{\frac{36}{5}}}
        {|u|^{\frac{36}{5}}+|v|^{\frac{36}{5}}},
                   \;\;\;              \text { if }\;\;|u|\geq 1, |v|\geq 1,\\
            40\left(\sin(2\pi x_1)+|\sin(2\pi x_1)|\right)^{\frac{11}{5}}
   \limsup_{|(u,v)|\rightarrow \infty}
   \frac{|u|^{6}|u|^{\frac{6}{5}}+|v|^{6}|v|^{\frac{6}{5}}+|u|^{\frac{36}{5}}+|v|^{\frac{36}{5}}}
        {|u|^{\frac{41}{5}}+|v|^{\frac{36}{5}}},
                   \;\;\;           \text { if }\;\;0\leq |u|< 1, |v|\geq 1,\\
           40\left(\sin(2\pi x_1)+|\sin(2\pi x_1)|\right)^{\frac{11}{5}}
   \limsup_{|(u,v)|\rightarrow \infty}
   \frac{|u|^{6}|u|^{\frac{6}{5}}+|v|^{6}|v|^{\frac{6}{5}}+|u|^{\frac{36}{5}}+|v|^{\frac{36}{5}}}
        {|u|^{\frac{36}{5}}+|v|^{\frac{41}{5}}},
                   \;\;\;             \text { if }\;\;|u|\geq 1, 0\leq |v|< 1
      \end{cases}
 \\
 & <  +\infty.
 \end{align*}}
So the condition $(F_4)$ holds. Then by Theorem 3.1, system (\ref{eq6}) has at least one nontrivial weak solution. If we let
\begin{eqnarray}\label{aaa1}
F(x,t,s)\equiv  |t|^{5}\ln(1+|t|)+|s|^{5}\ln(1+|s|)+|t|^{3}|s|^{3} \mbox{ for all }  x\in\mathbb R^N \mbox{ and } (t,s)\in\mathbb R^2,
\end{eqnarray}
 then by Theorem 3.9, system (\ref{eq6}) has infinitely many nontrivial weak solutions of high energy.

\section{Remark on the semi-trivial solutions of  \eqref{eq1}}\label{section 6}
 In Theorem 3.1 and Theorem 3.9,  we do not exclude the possibility of semi-nontrivial solutions. Hence, it is possible that the solutions of system are $(u,v)=(0,v)$ or $(u,v)=(u,0)$ which are called as  semi-nontrivial solutions. In general, it is not a simple work to rule out the semi-nontrivial solutions and some extra assumptions have to be added. We refer readers to \cite{Chang1}, \cite{Chang2} and \cite{Wu2016} for  related work. If we make the extra assumption $F(x,u,v)=F(x,|u|,|v|)$, by using Corollary 2.6 in \cite{Wu2016} and combing with the proofs of Theorem 3.9, it is easy to exclude the semi-trivial solutions, that is, system  \eqref{eq1} possesses infinitely many  non semi-trivial solutions.
 Especially, we can obtain that system (\ref{eq6}) with $F$ satisfying (\ref{aaa1}) has infinitely many  non semi-trivial solutions.

 \vskip2mm
\noindent{\bf Acknowledgements}
\par
This project is supported by Yunnan Ten Thousand Talents Plan Young and Elite Talents Project and Candidate Talents Training Fund of Yunnan Province, China (No: 2017HB016).


 \end{document}